\documentclass[12pt]{amsart}

\usepackage{amsmath,amsfonts,amsthm,amsopn,amssymb}
\usepackage{verbatim,wasysym,mathrsfs}
\usepackage[left=2.6cm,right=2.6cm,top=3.1cm,bottom=3.1cm]{geometry}
\usepackage{cite,marginnote}
\usepackage{color,enumitem,graphicx}
\usepackage[colorlinks=true,urlcolor=blue, citecolor=red,linkcolor=blue,
linktocpage,pdfpagelabels, bookmarksnumbered,bookmarksopen]{hyperref}
\usepackage[english]{babel}
\usepackage[T1]{fontenc}
\usepackage{lmodern}
\usepackage[hyperpageref]{backref}

\numberwithin{equation}{section}

\makeindex
\newtheorem{theorem}{Theorem}[section]
\newtheorem{proposition}[theorem]{Proposition}
\newtheorem{lemma}[theorem]{Lemma}
\newtheorem{remark}[theorem]{Remark}

\newtheorem{corollary}[theorem]{Corollary}
\newtheorem{definition}[theorem]{Definition}

\newcommand{\al}{\alpha}

\def\e{{\rm e}}

\newcommand{\lab}{\label}
\newcommand{\bt}{\begin{theorem}}
\newcommand{\et}{\end{theorem}}
\newcommand{\bl}{\begin{lemma}}
\newcommand{\el}{\end{lemma}}
\newcommand{\bd}{\begin{definition}}
\newcommand{\ed}{\end{definition}}
\newcommand{\bc}{\begin{corollary}}
\newcommand{\ec}{\end{corollary}}
\newcommand{\bp}{\begin{proof}}
\newcommand{\ep}{\end{proof}}
\newcommand{\bo}{\begin{proposition}}
\newcommand{\eo}{\end{proposition}}
\newcommand{\br}{\begin{remark}}
\newcommand{\er}{\end{remark}}
\newcommand{\be}{\begin{equation}}
\newcommand{\ee}{\end{equation}}
\newcommand{\ba}{\begin{align}}
\newcommand{\ea}{\end{align}}
\newcommand{\bn}{\begin{enumerate}}
\newcommand{\en}{\end{enumerate}}
\newcommand{\bg}{\begin{align*}}
\newcommand{\bcs}{\begin{cases}}
\newcommand{\ecs}{\end{cases}}

\newcommand{\Tail}{\text{Tail}}
\newcommand{\R}{\mathbb R}

\def\N{\mathbb N}

\makeatletter
\def\@makefnmark{}
\makeatother
\newcommand{\tC}{\widetilde C}
\newcommand{\tv}{\widetilde v}
\newcommand{\dd}{{\rm d}}
\newcommand{\bw}{\bar w}
\newcommand{\bean}{\begin{eqnarray*}}
\newcommand{\eean}{\end{eqnarray*}}
\newcommand{\loc}{\operatorname{\rm loc}}

\renewcommand{\epsilon}{\varepsilon}

\title[Nonlocal planar Schr\"odinger-Poisson systems]{Nonlocal planar Schr\"odinger-Poisson systems\\ in the fractional Sobolev limiting case}

\author[D. Cassani]{Daniele Cassani}
\author[Z. Liu]{Zhisu Liu}
\author[G. Romani]{Giulio Romani}

\address[Daniele Cassani and Giulio Romani]{
	\newline\indent Dipartimento di Scienza e Alta Tecnologia,
	\newline\indent Universit\`{a} degli Studi dell'Insubria
	\newline\indent and
	\newline\indent RISM-Riemann International School of Mathematics
	\newline\indent Villa Toeplitz, Via G.B. Vico, 46 - 21100 Varese, Italy}
	\email{\href{mailto:daniele.cassani@uninsubria.it}{daniele.cassani@uninsubria.it}}
\email{\href{mailto:giulio.romani@uninsubria.it}{giulio.romani@uninsubria.it}}

\address[Zhisu Liu]{\newline\indent Center for Mathematical Sciences/School of Mathematics and Physics,
\newline\indent
China University of Geosciences,
\newline\indent
Wuhan, Hubei, 430074, PR China
\newline\indent and
\newline\indent Dipartimento di Scienza e Alta Tecnologia,
\newline\indent Universit\`{a} degli Studi dell'Insubria
\newline\indent Villa Toeplitz, Via G.B. Vico, 46 - 21100 Varese, Italy}
\email{\href{mailto:liuzhisu@cug.edu.cn}{liuzhisu@cug.edu.cn}}

\thanks{Corresponding author: Daniele Cassani ({\tt daniele.cassani@uninsubria.it.})}

\thanks{Zhisu Liu is supported by the National Natural Science Foundation of China (No.~12226331), and the Fundamental Research Funds for the Central Universities, China University of Geosciences (Wuhan, Grant number: CUG2106211; CUGST2).}

\date{\today}
\subjclass[2020]{Primary: 35A15. Secondary: 35R11, 35J62, 35Q55, 35B09, 35B06.}
\keywords{Schr\"odinger-Newton system, limiting fractional Sobolev embeddings, Choquard type equations, p-fractional Laplacian, exponential growth, variational methods, positive solutions, moving planes and symmetry.}

\begin{document}

\begin{abstract} We study the nonlinear Schr\"odinger equation for the $s-$fractional $p-$Laplacian strongly coupled with the Poisson equation in dimension two and with $p=\frac2s$, which is the limiting case for the embedding of the fractional Sobolev space $W^{s,p}(\mathbb{R}^2)$. We prove existence of solutions by means of a variational approximating procedure for an auxiliary Choquard equation in which the uniformly approximated sign-changing logarithmic kernel competes with the exponential nonlinearity. Qualitative properties of solutions such as symmetry and decay are also established by exploiting a suitable moving planes technique. 
\end{abstract}

\maketitle

\section{Introduction}\label{sec1}

\noindent In this paper we investigate existence and symmetry of positive solutions for the following strong coupling of the $(s,\frac2s)-$fractional Schr\"odinger equation with the Poisson equation in the whole plane 
\begin{equation*}\tag{SP$_s$}\label{SPs}
	\begin{cases}
		(-\Delta)_\frac2s^su+|u|^{\frac2s-2}u=\phi\,f(u)\\
		-\Delta\phi=F(u)
	\end{cases}
\end{equation*}
where $s\in(0,1)$, $f$ is a nonnegative real function and $F$ its primitive vanishing at zero. The nonlocal operator in the Schr\"odinger equation is the $s-$fractional $p-$Laplacian with $p=\frac2s$, which is the so-called limiting case for the Sobolev embedding of the fractional Sobolev space (see Section \ref{Sec_2}) $$W^{s,p}(\R^2)\hookrightarrow L^{\frac{2p}{2-sp}}(\R^2)$$ 
and it is well-known that this implies the nonlinearity may grow exponentially at infinity.

\medskip

\noindent In the local case $s=1$, (SP$_s$) reduces to the deeply studied Schr\"odinger-Poisson system
\begin{equation*}\tag{SP}\label{SP}
	\begin{cases}
		-\Delta u+u=\phi\,f(u)\\
		-\Delta\phi=F(u)
	\end{cases}
\end{equation*}
which emerges in several fields of Physics: in electrostatics, it models the interaction of two identically charged particles; in quantum mechanics yields a model for the self-interaction of the wave function with its own gravitational field; it is also related with the Hartree model for crystals, see \cite{BF} and references therein. In the higher-dimensional case $N\geq3$, and when $f$ has a polynomial growth, there is an extensive literature about \eqref{SP}, see the survey \cite{MV} and references therein. A strategy to investigate \eqref{SP} is to move the attention to the corresponding Choquard equation
\begin{equation*}\tag{Ch}\label{Ch}
	-\Delta u+u=\left(I_N\ast F(u)\right)f(u)\quad\ \mbox{in}\ \R^N\,,
\end{equation*}
obtained by solving the Poisson equation in \eqref{SP} by means of convolution with the (positive) Riesz kernel
$$I_N:=\frac{C_N}{|x|^{N-2}},\qquad N\geq3$$
where $C_N$ is an explicit positive constant depending on the dimension $N$, and then formally (see the discussion carried out in Section \ref{Sec_equiv}) inserting $\phi=\phi_u:=I_N\ast F(u)$ in the Schr\"odinger equation. Besides the effect of variables reduction, the advantage of this approach is that the Choquard equation \eqref{Ch} can be studied by variational techniques. Moreover, new interesting phenomena arise, such as the appearance of a lower-critical exponent in addition to the usual upper-critical exponent, as in the Sobolev case, see \cite{CZ,CVZ} and references therein.

\noindent In dimension $N=2$ just a few results are available. In this case the Riesz kernel is logarithmic
\begin{equation}\label{Riesz_2}
	I_2:=\frac1{2\pi}\ln\frac1{|x|},
\end{equation}
and therefore unbounded both from below and from above, which together with the fact that it is sign changing introduces a major difficulty with respect to the higher dimensional case. In fact, the functional associated to the corresponding Choquard equation \eqref{Ch} may be not well-defined in the natural space $H^1(\R^2)$. When $f(u)=u$, the approach developed in \cite{CW,DW,BCV,CW2}, originating from the unpublished work of Stubbe \cite{Stubbe}, allows one to study the equation in a constrained subspace of $H^1(\R^2)$, in which the log-convolution term is well-defined. However, in dimension two it is well-known that $H^1(\R^2)\hookrightarrow L^q(\Omega)$ for all $q\geq1$ but not in $L^\infty(\Omega)$ and that  the maximal degree of summability for functions with membership in $H^1(\R^2)$ is exponential, namely $\int (e^{\alpha u^2}-1)\,\dd x$ stays bounded for all $\alpha>0$, see \cite{Ruf}. Several extensions and refinements to this result have been proposed, among which the extension to the higher-order Sobolev spaces $W^{m,p}$ with $m\in\N$, in the limiting case $N=mp$, see \cite{RS}.

\noindent Exponential nonlinearities in \eqref{Ch} when $N=2$, still mantaining the polynomial Riesz kernel and so loosing the connection to the Schr\"odinger-Poisson system, have been first considered in \cite{ACTY}. On the other hand, a special Schr\"odinger-Poisson system with logarithmic kernel and exponential nonlinearity, but not in gradient form as \eqref{SP}, has been studied in \cite{AF}. The tuning of the two phenomena: logarithmic growth of the Riesz kernel and maximal exponential growth  for system \eqref{SP} was first tackled in \cite{CT}. Here, the authors establish a proper functional setting by means of a log-weighted version of the Poho\v zaev-Trudinger inequality, so that the functional associated to \eqref{Ch} turns out to be well-defined. An extension of these techniques to quasilinear Schr\"odinger-Poisson systems was established in \cite{BCT}. A different approach has been recently proposed in \cite{LRZ, LRTZ,CDL}. Instead of working directly on the logarithmic Choquard equation, the authors consider a family of approximating problems, each of them involving only a polynomial Riesz kernel, and prove that the limit of a sequence of approximating solutions converges to a solution to the original problem. This approach has the advantage of working in the $H^1(\R^2)$ space context and one finds a posteriori that the logarithmic convolution term turns out to be well-defined at least on the solution. Let us finally mention that this method has been recently used in \cite{Ro} to cover also the zero-mass case. 

\noindent As a consequence of \cite{Parini19,Zhang19}, it is meaningful to study system \eqref{SPs} in presence of exponential nonlinearities. Fractional Schr\"odinger-Poisson systems and fractional Choquard equations have been recently studied in the subcritical regime, precisely in higher dimension $N>sp$ and with polynomial nonlinearities, see \cite{dASS, LM, HW, GRW}. To the very best of our knowledge, the only results for fractional Choquard equations with exponential nonlinearities are obtained in \cite{CDB} in the one dimensional case, in \cite{BM} where the logarithmic kernel and the exponential nonlinearity are not combined, and in \cite{YRTZ} where the Riesz kernel is polynomial; see also \cite{YRCW} for related results. 

\medskip

\noindent Here we address both the problem of proving existence and establishing qualitative properties of solutions to the fractional planar system \eqref{SPs}, in the limiting case of logarithmic kernel and exponential nonlinearity. The first nontrivial step is proving that solutions of the system can be obtained by solving the related Choquard equation
\begin{equation}\label{Chs}\tag{Ch$_s$}
	(-\Delta)^s_{\frac2s} u+|u|^{\frac2s-2}u=\frac1{2\pi}\left(\ln\frac1{|x|}\ast F(u)\right)f(u)\quad\text{in}\ \,\R^2\,.
\end{equation}

\medskip

\noindent Before stating our main results, let us briefly introduce the functional setting and make precise the definition of solution we deal with for \eqref{SPs} and \eqref{Chs}. \noindent The $(s,\tfrac2s)-$fractional Laplace operator is pointwisely defined as
$$
(-\Delta)^s_{\frac2s}u(x):=2\,PV\!\int_{\R^2}\frac{|u(x)-u(y)|^{\frac2s-2}(u(x)-u(y))}{|x-y|^4}\,{\rm d}y,\quad x\in\R^2,
$$
where \textit{PV} stands for the Cauchy Principal Value which is well-defined for all $x\in\R^2$, for functions in $C^{1,1}_{loc}(\R^2)$ which enjoy suitable integrability conditions at infinity, see Section \ref{Sec_equiv} and \cite[Lemma 5.2]{Chen18}. 

\noindent The fractional Sobolev space $W^{s,\frac2s}(\R^2)$ is defined by
$$
W^{s,\frac2s}(\R^2):=\{u\in L^{\frac2s}(\R^2): [u]_{s,\frac2s}<+\infty\},
$$
where $[u]_{s,\frac2s}$ denotes the Gagliardo seminorm, which in this context reads as follows 
$$
[u]_{s,\frac2s}=\bigg(\int_{\R^2}\int_{\R^2}\frac{|u(x)-u(y)|^{\frac2s}}{|x-y|^4}\,{\rm d}x\,{\rm d}y\bigg)^{\frac s2}.
$$
The fractional Sobolev-Slobodeckij space $W^{s,\frac2s}(\R^2)$ is a uniformly convex Banach space with norm 
$$
\|u\|:=\bigg(\|u\|_{\frac2s}^{\frac2s}+[u]_{s,\frac2s}^{\frac2s}\bigg)^{\frac s2}.
$$

\noindent For $\gamma>0$ the weighted Lebesgue space $L_\gamma(\R^2)$ is defined as
$$L_\gamma(\R^2):=\Big\{u\in L^1_{\loc}(\R^2)\,\Big|\,\int_{\R^2}\frac{|u(x)|}{1+|x|^{2+2\gamma}}\,{\rm d}x<+\infty\Big\}\,.$$

\begin{definition}\label{sol_Poisson}
	For $f\in\mathcal S'(\R^2)$ we say that a function $\phi\in L_1(\R^2)$ is a solution of the linear Poisson equation $-\Delta\phi=f$ in $\R^2$ if
	$$\int_{\R^2}\phi(-\Delta\varphi)=\langle f,\varphi\rangle\qquad\mbox{for all}\ \,\varphi\in\mathcal S(\R^2)\,.$$
\end{definition}

\begin{definition}[Solution of \eqref{SPs}]\label{sol_SP}
	We say that $(u,\phi)$ is a weak solution of the Schr\"odinger-Poisson system \eqref{SPs} if
	\begin{equation*}
		\int_{\R^2}\int_{\R^2}\frac{|u(x)-u(y)|^{\frac2s-2}(u(x)-u(y))(\varphi(x)-\varphi(y))}{|x-y|^4}\,{\rm d}x\,{\rm d}y+\int_{\R^2}|u|^{\frac2s-2}u\varphi\,{\rm d}x=\int_{\R^2}\phi f(u)\varphi\,{\rm d}x
	\end{equation*}
	for all $\varphi\in W^{s,\frac2s}(\R^2)$, and $\phi$ solves $-\Delta\phi=F(u)$ in $\R^2$ in the sense of Definition \ref{sol_Poisson}.
\end{definition}

\begin{definition}[Solution of \eqref{Chs}]\label{sol_Choquard_log}
	We say that $u\in W^{s,\frac2s}(\R^2)$ is a \textit{weak solution of} \eqref{Chs} if
	\begin{equation}\label{sol_Choquard_log_test}
		\begin{split}
			\int_{\R^2}\int_{\R^2}&\frac{|u(x)-u(y)|^{\frac2s-2}(u(x)-u(y))(\varphi(x)-\varphi(y))}{|x-y|^4}\,{\rm d}x\,{\rm d}y+\int_{\R^2}|u|^{\frac2s-2}u\varphi\,{\rm d}x\\
			&=\frac1{2\pi}\int_{\R^2}\left(\int_{\R^2}\ln\frac1{|x-y|}F(u(y))\,{\rm d}y\right)f(u(x))\varphi(x)\,{\rm d}x
		\end{split}
	\end{equation}
	for all $\varphi\in W^{s,\frac2s}(\R^2)$.
\end{definition}

\paragraph{\textbf{Assumptions.}}
Throughout the paper, we consider the following assumptions on the nonlinearity $f$:
\begin{itemize}
\item[($f_1$)] $f\in C^1(\R)$, $f\geq 0$, $f(t)=0$ for $t\leq0$, and $f(t)=o(t^{\frac2s-1})$ as $t\to0^+$;
\item[($f_2$)] there exist constants $b_1,b_2>0$
such that for any $t>0$,
$$0<f(t)\leq b_1+b_2\Phi_{2,s}(\alpha_*|t|^{\frac{2}{2-s}})\,,$$
where $\Phi_{2,s}(\,\cdot\,)$ has exponential growth, see Proposition \ref{Thm:m-theo}; 
\item[($f_3$)] there exists $\tau\in(0,s)$ such that
	$$1-s+\tau\leq\frac{F(t)f'(t)}{f^2(t)}<1+\mu(s,\tau)\quad\mbox{for any}\ \,t>0,$$
	where $\mu(s,\tau)$  is explicitly given in Lemma \ref{Lem:integral-F};
\item[($f_4$)] $\lim\limits_{t\to+\infty}\frac{F(t)f'(t)}{f^2(t)}=1$ or equivalently
$\lim\limits_{t\to+\infty}\frac{\rm d}{{\rm d}t}\frac{F(t)}{f(t)}=0\,$;
\item[($f_5$)] there exists $\beta_0>1$ depending on $s$ such that
	\begin{equation*}
		f(t)F(t)\geq\beta\,t^\frac2s\quad\ \mbox{for all}\ \ t>T(s),
	\end{equation*}
	for some $\beta>\beta_0$. The values of $\beta_0$ and $T(\,\cdot\,)$ are explicitly given in Lemma \ref{Lem:guji}.  
\end{itemize}

\medskip

\noindent \textbf{Main results.} We are now in the position to state our main results. The first one, proved in Section \ref{Sec_equiv}, concerns the relationship between the solutions of the Choquard equation \eqref{Chs} and the corresponding Schr\"odinger-Poisson system \eqref{SPs}: this is the fractional counterpart of \cite[Theorem 2.1]{BCT}. In particular we make rigorous the fact that from a weak solution $u$ of \eqref{Chs} one obtains a weak solution to \eqref{SPs}. 

\begin{theorem}[\eqref{Chs}$\,\,\Longrightarrow\,\,$\eqref{SPs}]\label{thm_equiv}
	Suppose the nonlinearity $f$ satisfies ($f_1$)-($f_2$). Let $u\in W^{s,\frac2s}(\R^2)$ be a positive weak solution of the Choquard equation \eqref{Chs} in the sense of Definition \ref{sol_Choquard_log} and define $\phi_u:=I_2\ast F(u)$. Then $(u,\phi_u)\in W^{s,\frac2s}(\R^2)\times L_{\frac12}(\R^2)$ is a solution of the Schr\"odinger-Poisson system \eqref{SPs} in the sense of Definition \ref{sol_SP}. Moreover $u\in L^\infty(\R^2)\cap C^{0,\nu}_{\loc}(\R^2)$ for some $\nu\in(0,1)$, and
	\begin{equation}\label{asympt_phi}
		\phi_u(x)=-\frac{\|F(u)\|_1}{2\pi}\ln|x|+o(1)\qquad\mbox{as}\ \ |x|\to+\infty\,.
	\end{equation}
\end{theorem}

\noindent 	Note that we do not prove that the two problems are equivalent, namely that they have the same set of solutions. Finding a proper functional setting in which this holds is still open, even for the local case \eqref{SP}--\eqref{Ch}, see \cite[Section 2]{BCT}.

\noindent The system \eqref{SPs} is autonomous, thus we expect that positive solutions enjoy radial symmetry. In the next result, by exploiting the connection between \eqref{Chs} and \eqref{SPs} established in Theorem \ref{thm_equiv}, we show that this is indeed the case for solutions which are regular enough. The regularity will be needed to have a pointwisely defined $(s,\frac2s)$-fractional Laplace operator and to be able to exploit a moving plane technique. This method has been firstly employed to $(s,p)$-fractional operators by Chen and Li in \cite{Chen18} for the single Schr\"odinger equation; here, we further develop it to cover the case of systems. 

\begin{theorem}[Symmetry for \eqref{Chs}]\label{Thm:symmetric}
	Suppose ($f_1$)-($f_2$) are satisfied, and let $u\in C_{loc}^{1,1}(\R^2)\cap W^{s,\frac2s}(\R^2)$ be a positive solution of \eqref{Chs}. Then $u$ is radially symmetric about the origin and monotone decreasing.
\end{theorem}

\noindent Observe that in Theorems \ref{thm_equiv} and \ref{Thm:symmetric} we only need that $f$ is subcritical or critical in the sense of Proposition \ref{Thm:m-theo}. In order to establish the existence of solutions for \eqref{Chs} -- and therefore for \eqref{SPs} by Theorem \ref{thm_equiv} -- we require the full set of assumptions ($f_1$)--($f_5$). 

\begin{theorem}[Existence for \eqref{Chs}]\label{Thm:duozhjie}
	Suppose $f$ satisfies ($f_1$)-($f_5$), then \eqref{Chs} possesses a positive radially symmetric solution $u\in W^{s,\frac2s}(\R^2)$ such that
	\begin{equation}\label{logFF}
		\bigg|\int_{\R^2}\bigg(\ln\frac1{|x|}\ast F(u)\bigg)F(u)\,{\rm d}x\bigg|<+\infty\,.
	\end{equation}
\end{theorem}

\noindent As a consequence of Theorems \ref{thm_equiv} and \ref{Thm:duozhjie}, we have 
\begin{corollary}[Existence for \eqref{SPs}]
	Suppose $f$ satisfies ($f_1$)-($f_5$), then \eqref{SPs} possesses a solution $(u,\phi)\in W^{s,\frac2s}(\R^2)\times L_\frac12(\R^2)$ such that:
	\begin{itemize}
		\item[i)] $u\in L^\infty(\R^2)\cap C^{0,\nu}_{\loc}(\R^2)$ for some $\nu\in(0,1)$, is positive, radially symmetric and \eqref{logFF} holds;
		\item[ii)] $\phi=\phi_u:=I_2\ast F(u)$ and the asymptotic behaviour \eqref{asympt_phi} is satisfied.
	\end{itemize}
\end{corollary}

\vskip0.2truecm

\noindent \textbf{Overview.} In the next section we give motivations and discuss consequences of our assumptions together with some preliminaries. In Section \ref{Sec_equiv} we study the relationship between solutions of the Choquard equation \eqref{Chs} and the related Schr\"odinger-Poisson system \eqref{SPs} which is quite delicate as it depends on the notion of solution and regularity issues. Apparently there is no equivalence in general, even in the local case, see \cite{BCT}; Theorem \ref{thm_equiv} is a step forward towards a complete understanding of this phenomenon. In Section \ref{Sec_MVP} we prove the symmetry result of Theorem \ref{Thm:symmetric} by extending the moving-plane technique of \cite{Chen18}. Finally, in Section \ref{Sec_existence} we exploit all previous results to prove Theorem \ref{Thm:duozhjie} by means of a variational approach together with a uniform asymptotic approximation technique. 

\medskip

\paragraph{\textbf{Notation.}} For $R>0$ and $x_0\in\R^2$ we denote by $B_R(x_0)$ the ball of radius $R$ and center $x_0$. Given a set $\Omega\subset\R^2$, we denote $\Omega^c:=\R^2\setminus\Omega$, and its characteristic function as $\chi_\Omega$. The space of the infinitely differentiable functions which are compactly supported is denoted by $C^\infty_0(\R^2)$, while $L^p(\R^2)$ with $p\in[1,+\infty]$ is the Lebesgue space of $p$-integrable functions. The norm of $L^p(\R^2)$ is denoted by $\|\cdot\|_p$. The spaces $C^{0,\nu}(\R^2)$ for $\nu\in(0,1)$ are usual spaces of H\"older continuous functions. The space $\mathcal{S}$ is the Schwartz space of rapidly decreasing functions and $\mathcal{S}'$ the dual space of tempered distributions. For $q>0$ we define $q!:=q(q-1)\cdots(q-\lfloor q\rfloor)$, where $\lfloor q\rfloor$ denotes the largest integer strictly less than $q$; if $q>1$ its conjugate H\"older exponent is $q':=\frac q{q-1}$. Finally, $o_n(1)$ denotes a vanishing real sequence as $n\to+\infty$. Hereafter, the letter $C$ will be used to denote positive constants which are independent of relevant quantities and whose value may change from line to line.

\section{Preliminaries}\setcounter{equation}{0}\label{Sec_2}
\noindent The Poho\v zaev-Trudinger inequality has been extended to the Sobolev fractional setting by Parini and Ruf \cite{Parini19} for bounded domains and by Zhang \cite{Zhang19} for the whole space, results which we resume below. Let    
$$\Phi_{2,s}(t):=\e^t-\sum_{j=0}^{j_{\frac2s}-2}\frac{t^j}{j!}\,,$$
for $t\geq 0$, where $j_{\frac2s}:=\min\{j\in\N:\,j\geq\frac2s\}$.
\begin{proposition}[\cite{Parini19, Zhang19} ]\label{Thm:m-theo}
	Let $s\in(0,1)$, then for all $\alpha>0$ one has
	\begin{equation}\label{Trudinger}
		\int_{\R^2}\Phi_{2,s}(\alpha|u|^{\frac2{2-s}})\,{\rm d}x<+\infty\,.
	\end{equation}
	Moreover, there exists $\alpha_*\in(0,+\infty)$ such that for $0\leq\alpha<\alpha_*$,
	$$
	\sup\limits_{u\in W^{s,\frac2s}(\R^2),\,\|u\|_{W^{s,\frac2s}(\R^2)}\leq1}\int_{\R^2}\Phi_{2,s}(\alpha|u|^{\frac2{2-s}})\,{\rm d}x<+\infty\,.
	$$
	Hence
	$$
	\alpha_*:=\sup\bigg\{\alpha:\,\sup\limits_{u\in W^{s,\frac2s}(\R^2),\,\|u\|_{W^{s,\frac2s}(\R^2)}\leq1}
	\int_{\R^2}\Phi_{2,s}(\alpha|u|^{\frac2{2-s}})\,{\rm d}x<+\infty\bigg\}\,.
	$$
	Moreover, for $\alpha>\alpha_{s,2}^*\,$,
	$$
	\sup\limits_{u\in W^{s,\frac2s}(\R^2),\,\|u\|_{W^{s,\frac2s}(\R^2)}\leq1}\int_{\R^2}\Phi_{2,s}(\alpha|u|^{\frac2{2-s}})\,{\rm d}x=+\infty\,,
	$$
	where
	$$
	\alpha_{s,2}^*:=2^\frac{4+s}{2-s}\bigg(\pi^2\Gamma(1+\tfrac2s)\sum_{k=0}^{+\infty}\frac1{(1+k)^{2/s-1}}\bigg)^{\frac s{2-s}}.
	$$
\end{proposition}

\begin{remark}As remarked in \cite{Parini19,Zhang19}, the result of Proposition \ref{Thm:m-theo} is not sharp in the sense of Moser, since $\alpha_{s,2}^*$ is just an upper bound for the sharp exponent $\alpha_*$, while obtaining the precise value of $\alpha_*$ is still an open problem.
\end{remark}
\noindent Let us point out some immediate consequences of the assumptions ($f_1$)--($f_5$):

\begin{remark}\label{Rmk_ass}
\begin{enumerate}
\item[]
	\item[(i)] From ($f_1$) and ($f_2$) it is easy to infer that the exists a constant $C$ such that
	$$F(t)\leq C\left(t^\frac2s+t^\frac2s\Phi_{2,s}(\alpha_*|t|^{\frac2{2-s}})\right)\quad\mbox{for all}\ \,t>0\,;$$
	\item[(ii)] Assumption $(f_3)$ implies that $f$ is monotone increasing. Moreover,
	$$\frac{{\rm d}}{{\rm d}t}\frac{F(t)}{f(t)}=\frac{f^2(t)-F(t)f'(t)}{f^2(t)}\leq s-\tau\,,$$
	from which one infers
	$$F(t)\leq(s-\tau)tf(t)\quad\ \mbox{for any}\ \ t\geq0\,,$$
	a Ambrosetti-Rabinowitz type condition. The upper bound in ($f_3$) is fulfilled by nonlinearities of the kind $t\mapsto e^{t^\alpha}$ with $\alpha>1$, in particular by the critical growth $t\mapsto e^{t^\frac2{2-s}}$;
	\item[(iii)] Assumption ($f_4$) requires that the function $f$ grows exponentially at infinity. Moreover, from ($f_4$) one may deduce that for any $\varepsilon>0$ there exists $M_\varepsilon$ such that for $t\geq M_\varepsilon$
	\begin{equation}\label{Rmk_ass_2_eq}
		F(t)\leq\varepsilon tf(t)\,.
	\end{equation}
	Indeed, there exists $N_\varepsilon>0$ such that $\frac{{\rm d}}{{\rm d}s}\frac{F(s)}{f(s)}<\frac\varepsilon2$ for all $s>N_\varepsilon$, and integrating over $[N_\varepsilon,t]$ one gets
	$$\frac{F(t)}{f(t)}-\frac{F(N_\varepsilon)}{f(N_\varepsilon)}=\int_{N_\varepsilon}^t\frac{{\rm d}}{{\rm d}s}\frac{F(s)}{f(s)}\,{\rm d}s<\frac\varepsilon2(t-N_\varepsilon)\,.$$
	Then \eqref{Rmk_ass_2_eq} follows, by restricting to $t>M_\varepsilon\geq N_\varepsilon$, where $M_\varepsilon$ is chosen such that $\frac{F(N_\varepsilon)}{f(N_\varepsilon)}-\varepsilon N_\varepsilon<\frac\varepsilon2M_\varepsilon\,$;
	\item[(iv)] ($f_5$) is a condition about middle-range values of $f(t)$, since by ($f_4$) it is automatically verified at infinity. It can be related to the De Figueiredo-Miyagaki-Ruf condition \cite{Figueiredo95}, and is crucial in order to estimate the mountain pass level and gain compactness, see Lemma \ref{Lem:guji}. We point out that the fractional analogue of the Moser sequence, usually exploited to estimate the mountain pass level, is here spoiled by the fact that it is not known whether the value $\alpha_*$ coincides with $\alpha_{s,2}$ in Proposition \ref{Thm:m-theo}. Assumption ($f_5$), combined with the choice of a simple fixed test function, which allows explicit computations of its seminorm, will enable us to give an explicit estimate of the mountain pass level.
	\item[(v)] Finally let us give an example of nonlinearity satisfying ($f_1$)--($f_5$). Let $A,B,t_*>0$ and define
\begin{equation*}
	F(t)=\begin{cases}
		A\,t^{\frac2s+1}\quad&\mbox{for}\ \,t\in(0,t_*),\\%
		B\,\e^{t^\frac2{2-s}}\quad&\mbox{for}\ \,t\geq t_*,
	\end{cases}
\end{equation*}
and extend it by $0$ on $\R^-$. Note that $F$ has critical growth at $\infty$, and that assumptions ($f_1$) and ($f_2$) are readily satisfied for any choice of the constants. Moreover, for $t\in(0,t_*)$ one has
$$\frac{F(t)f'(t)}{f^2(t)}=1-\frac s{2+s},$$
which is bigger than or equal to $1-s+\tau$, e.g. with the choice $\tau=\tau(s):=\frac{s(1+s)}{2+s}$. Hence the lower bound in ($f_3$) is satisfied on $(0,t_*)$. Assumption ($f_4$) is also easily verified, since
\begin{equation*}
	\lim_{t\to+\infty}\frac{F(t)f'(t)}{f^2(t)}=\lim_{t\to+\infty}\left(1+\frac{s}{2}t^{-\frac{2}{2-s}}\right)=1\,.
\end{equation*}
By the same computation we also note that the lower bound in ($f_3$) is trivially satisfied on $(t_*,+\infty)$, whereas the upper bound in ($f_3$) holds, provided 
$$\frac{s}{2}t^{-\frac{2}{2-s}}<\mu(s,\tau(s))\,,$$
with $\mu(s,\tau)$ given in Lemma \ref{Lem:integral-F}, that is, for $t>t_*$ where 
$$t_*:=\left(\frac{s}{2\mu(s,\tau(s))}\right)^\frac{2-s}2.$$ 
\noindent In order for $f=F'$ to be continuous on $\R$, one needs to impose a linear dependence between the constants $A$ and $B$. Finally, ($f_5$) is also satisfied, provided one chooses $A$ and $B$ large enough. This is a working example, however, let us point out that these restrictions on constants can be relaxed. 
\end{enumerate}
\end{remark}

\medskip

\noindent For $s\in(0,1)$ and $\Omega\subset\R^2$ we define
$$W_0^{s,\frac2s}(\Omega):=\left\{u\in W^{s,\frac2s}(\R^2)\,|\,u_{|_{\R^2\setminus\Omega}}\equiv 0\right\}.$$
The subspace of radial functions in $W^{s,\frac2s}(\R^2)$ is defined as
$$W_r^{s,\frac2s}(\R^2):=\{u\in W^{s,\frac2s}(\R^2)\,|\,\,u(x)=u(|x|)\}\,.$$
The following compactness result due to Lions \cite{Lions82} will be crucial in our analysis.
\begin{lemma}\label{Lem:lions}
	Let $s\in(0,1]$. Then $W_r^{s,\frac2s}(\R^2)$ is compactly embedded in $L^q(\R^2)$ for any $\frac2s<q<\infty$.
\end{lemma}

\noindent Let us also recall the well-known Hardy-Littlewood-Sobolev inequality, see \cite[Theorem 4.3]{LL}, which will be frequently used throughout the paper.
\begin{lemma}(Hardy-Littlewood-Sobolev inequality)\label{Lem:HLS}
	Let $N\geq1$, $q,r>1$, and $\alpha\in(0,N)$ with $\tfrac1q+\tfrac\alpha N+\tfrac1r=2$. There exists a constant $C=C(N,\alpha,q,r)$ such that for all $f\in L^q(\R^N)$ and $h\in L^r(\R^N)$ one has
	$$
	\int_{\R^N}\left(\frac1{|x|^\alpha}\ast f(x)\right)h(x)\,{\rm d}x\leq C\|f\|_q\|h\|_r\,.
	$$
\end{lemma}

\noindent The following technical lemma will be useful in the sequel,
\begin{lemma}\label{Lem:jianjin}
	Assume $u\in W^{s,\frac2s}(\R^2)$ and ($f_1$)-($f_2$) hold. Then for any $0<\kappa\leq1$ we have
	$$
	\int_{\R^2}\frac{F(u(y))}{|x-y|^{\kappa}}\,{\rm d}y-\frac{1}{|x|^{\kappa}}\int_{\R^2}F(u(y))\,{\rm d}y\to0\,,\quad\text{as}\ \ |x|\to+\infty\,.
	$$
\end{lemma}
\begin{proof}
	First, we have $F(u)\in L^1(\R^2)$ by Remark \ref{Rmk_ass}(i) and Proposition \ref{Thm:m-theo}. Consider now $|x|\geq1$ and note that
	$$\int_{\R^2}\frac{F(u(y))}{|x-y|^{\kappa}}\,{\rm d}y-\frac1{|x|^{\kappa}}\int_{\R^2}F(u(y))\,{\rm d}y=\int_{\R^2}\bigg(\frac1{|x-y|^{\kappa}}-\frac1{|x|^{\kappa}}\bigg)F(u(y))\,{\rm d}y\,.$$
	For $x,y\in\R^2$ define $g(x,y):=\frac1{|x-y|^{\kappa}}-\frac1{|x|^{\kappa}}$. Observe that $g(x,y)\to0$ as $|x|\to+\infty$ for every $y\in\R^2$. Moreover,
	$$
	-1\leq g(x,y)\chi_{|x-y|\geq\frac12}(y)\leq2^{\kappa}\quad\text{for}\ \,x,y\in\R^2\ \,\text{with}\ \,|x|\geq1\,.
	$$
	Hence, Lebesgue's theorem implies that
	$$
	\int_{|y-x|\geq\frac12}g(x,y)F(u(y))\,{\rm d}y\to0\quad\text{as}\ \,|x|\to+\infty\,.
	$$
	Moreover, we have
	$$
	0\leq \frac{1}{|x|^{\kappa}}\int_{|x-y|\leq \frac12}F(u(y))\,{\rm d}y\leq\frac1{|x|^{\kappa}}\int_{\R^2}F(u(y))\,{\rm d}y\to0\quad\text{as}\ \,|x|\rightarrow+\infty\,,
	$$
	and
	\begin{equation}\lab{eqn:++}
		\aligned
			0&\leq\int_{|x-y|\leq\frac12}\frac1{|x-y|^{\kappa}}F(u(y))\,{\rm d}y\\
			&\leq\bigg(\int_{B_\frac12(x)}\frac1{|y|^{\kappa \sigma}}\,{\rm d}y\bigg)^{\frac1\sigma}\bigg(\int_{\R^2}F(u(y))^{\sigma'}\,{\rm d}y\bigg)^{\frac1{\sigma'}}\to0\quad\text{as}\ \,|x|\to+\infty\,,
		\endaligned
	\end{equation}
	by Theorem \ref{Thm:m-theo}, for $\sigma>1$ large enough and $\sigma'=\frac{\sigma}{\sigma-1}$. Based on the above estimates, we infer that
	$$
	\int_{\R^2}g(x,y)F(u(y))\,{\rm d}y\to0\quad\text{as}\ \,|x|\to\infty\,.
	$$
\end{proof}

\section{On the relationship between the nonlocal Choquard equation and the Schr\"odinger-Poisson system: proof of Theorem \ref{thm_equiv}}\label{Sec_equiv}

\noindent It is commonly used in the literature \cite{CW,CW2,CT,DW,AF} that solutions $u$ of the Choquard equation \eqref{Chs} correspond to solutions $(u,\phi)$ of the Schr\"odinger-Poisson system \eqref{SPs} by formally inverting the Laplace operator in the second equation and considering the only solution $\phi$ given by the convolution of the nonlinearity with the Riesz kernel, namely $\phi=\phi_u:=I_N\ast F(u)$. However, the Sobolev limiting case $N=sp$ is affected by regularity issues which have roots in the lack of regularity of the Fourier transform, see \cite{H}.  This turns out to be a quite delicate situation which we redeem in Theorem \ref{thm_equiv} where we settle a suitable functional framework in order to obtain solutions for \eqref{SPs} from solutions of the Choquard equation \eqref{Chs}. The complete equivalence between \eqref{SPs} and \eqref{Chs} is still open even in the non fractional case, see \cite{BCT}.

\noindent As a first step to prove Theorem \ref{thm_equiv}, we show that solutions to \eqref{Chs} are bounded by means of a Nash-Moser iteration argument in the spirit of \cite{AI,Bisci2022}. This fact, together with a suitable polynomial estimate of the decay at infinity, yields a log-weighted $L^1$-estimate for $F(u)$, which is the crucial ingredient to prove Theorem \ref{thm_equiv}. Then, we end up with $C^{0,\nu}_{\loc}(\R^2)$ for some $\nu\in(0,1)$, by standard regularity results. Throughout this section we assume ($f_1$)-($f_2$).

\begin{remark}\label{Rmk_positivity}
	Note that any nontrivial weak solution $u\in W^{s,\frac2s}(\R^2)$ of \eqref{Chs} is positive. Indeed, take $u^-:=\min\{u,0\}\in W^{s,\frac2s}(\R^2)$ as test function in \eqref{sol_Choquard_log_test}, to get
	\begin{equation*}
		\aligned
			0=&\int_{\R^2}\int_{\R^2}\frac{|u(x)-u(y)|^{\frac2s-2}(u(x)-u(y))(u^-(x)-u^-(y))}{|x-y|^4}\,{\rm d}x\,{\rm d}y+\int_{\R^2}|u|^{\frac2s-2}u u^-\,{\rm d}x\\
			&-\frac1{2\pi}\int_{\R^2}\left(\int_{\R^2}\ln\frac1{|x-y|}F(u(y))\,{\rm d}y\right)f(u(x))u^-(x)\,{\rm d}x\\
			\geq&\int_{\R^2}\int_{\R^2}\frac{|u^-(x)-u^-(y)|^\frac2s}{|x-y|^4}\,{\rm d}x\,{\rm d}y+\int_{\R^2}|u^-|^\frac2s\,{\rm d}x\,,
		\endaligned
	\end{equation*}
	which implies that $u^-\equiv0$. Hence, using the strong maximum principle for the $p-$fractional Laplacian \cite[Theorem 1.4]{Del17} we can infer that $u>0$ in $\R^2$.
\end{remark}

\begin{lemma}\label{Lemma_infty}
	Let $u\in W^{s,\frac2s}(\R^2)$ be a positive solution of the Choquard equation \eqref{Chs}, then $u\in L^\infty(\R^2)$.
\end{lemma}
\begin{proof}
	Since $f(u)\geq0$, for a nonnegative test function $\varphi\in W^{s,\frac2s}(\R^2)$ we have
	\begin{equation}\label{eqn:re1}
		\begin{split}
			\int_{\R^2}\int_{\R^2}&\frac{|u(x)-u(y)|^{\frac2s-2}(u(x)-u(y))(\varphi(x)-\varphi(y))}{|x-y|^4}\,{\rm d}x\,{\rm d}y+\int_{\R^2}|u|^{\frac2s-2}u\varphi\,{\rm d}x\\
			&=\frac1{2\pi}\int_{\R^2}\left(\int_{\R^2}\ln\frac1{|x-y|}F(u(y))\,{\rm d}y\right)f(u(x))\varphi(x)\,{\rm d}x\\
			&\leq\frac1{2\pi}\int_{\R^2}\left(\int_{\R^2}\frac{F(u(y))}{|x-y|}\,{\rm d}y\right)f(u(x))\varphi(x)\,{\rm d}x\\
			&\leq C\int_{\R^2}f(u(x))\varphi(x)\,{\rm d}x\,,
		\end{split}
	\end{equation}
	where in the last inequality we have used the fact that $\frac1{|\cdot|}\ast F(u)$ is bounded, which follows arguing as in \cite[Lemma 4.1]{ACTY}. For any $L>0$ and $\gamma>1$ define
	\begin{equation}\label{psi_Psi_K}
		\psi(t):=t(\min\{t,L\})^{\frac2s(\gamma-1)},\quad\ \Psi(t):=\int_{0}^t(\psi'(\tau))^{\frac s2}\,{\rm d}\tau,\quad\ K(t):=\frac s2|t|^{\frac2s}
	\end{equation}
	for $t>0$. By Jensen's inequality it is possible to show that
	\begin{equation}\label{eqn:re2-0}
		|\Psi(a)-\Psi(b)|^{\frac2s}\leq K'(a-b)(\psi(a)-\psi(b))
	\end{equation}
	holds for any $a,b\in\R^+$. Therefore
	\begin{equation*}
		\aligned
		|\Psi(u(x))-\Psi(u(y))|^{\frac2s}\leq|u(x)-u(y)|^{\frac2s-2}(u(x)-u(y))[u(x)u_L(x)^{\frac2s(\gamma-1)}-u(y)u_L(y)^{\frac2s(\gamma-1)}]\,,
		\endaligned
	\end{equation*}
	where $u_L:=\min\{u,L\}$. Taking $\psi(u)=u u_L^{\frac2s(\gamma-1)}$ as test function in \eqref{eqn:re1}, we obtain
	\begin{equation}\label{eqn:re3}
		\begin{split}
			[\Psi(u)]_{s,\frac2s}^{\frac2s}&\leq\int_{\R^2}\int_{\R^2}\frac{|u(x)-u(y)|^{\frac2s-2}(u(x)-u(y))[u(x)u_L(x)^{\frac2s(\gamma-1)}-u(y)u_L(y)^{\frac2s(\gamma-1)}]}{|x-y|^4}\,{\rm d}x\,{\rm d}y\\
			&\leq C\int_{\R^2}f(u)u u_L^{\frac2s(\gamma-1)}\,{\rm d}x-\int_{\R^2}|u|^{\frac2s}u_L^{\frac2s(\gamma-1)}\,{\rm d}x\,.
		\end{split}
	\end{equation}
	Recall that from ($f_1$)-($f_2$) we have that for any $\epsilon>0$, there exists $C_\epsilon>0$ such that
	$$
	f(u)\leq\epsilon |u|^{\frac2s-1}+C_{\epsilon}|u|^{\frac2s-1}\Phi_{2,s}(\alpha_*|u|^{\frac{2}{2-s}})\,.
	$$
	Hence, from \eqref{eqn:re3} we infer
	\begin{equation}\label{eqn:re4}
		\begin{split}
			{[\Psi(u)]}^{\frac2s}_{s,\frac2s}&\leq (C\epsilon-1)\int_{\R^2}|u|^{\frac2s} u_L^{\frac2s(\gamma-1)}{\rm d}x+C_{\epsilon}\int_{\R^2}\Phi_{2,s}(\alpha_*|u|^{\frac2{2-s}})|u|^{\frac2s}u_L^{\frac2s(\gamma-1)}{\rm d}x\\
			&\leq (C\epsilon-1)\int_{\R^2}|u|^{\frac2s} u_L^{\frac2s(\gamma-1)}{\rm d}x\\
			&\quad+C_{\epsilon}\bigg(\int_{\R^2}|u|^{\frac{2p}{s}}u_L^{\frac{2p}s(\gamma-1)}{\rm d}x\bigg)^{\frac1p}\!\bigg(\int_{\R^2}\Phi_{2,s}(\alpha_*p'|u|^{\frac2{2-s}})\,{\rm d}x\bigg)^{\frac1{p'}}\\
			&\leq (C\epsilon-1)\|u u_L^{\gamma-1}\|_{\frac2s}^{\frac2s}+C_{\epsilon}\|u u_L^{\gamma-1}\|_{\frac{2p}s}^{\frac2s}\,,
		\end{split}
	\end{equation}
	for $p>1$ and $p'=\frac p{p-1}$ by ($f_2$) and Theorem \ref{Thm:m-theo}. By \eqref{eqn:re2-0} one has $\Psi(u)\leq u u_L^{\gamma-1}$. Moreover,
	$$\Psi(u)\geq\frac1{\gamma}u u_L^{\gamma-1},$$
	due to the definition of $\Psi$. Indeed, on the one hand if $u(x)<L$ one has $\psi(u)(x)=u(x)^{\frac2s(\gamma-1)}\,$, thus 
	\begin{equation*}
		\Psi(u)(x)=\left(\frac2s(\gamma-1)+1\right)^{\frac s2}\int_0^{u(x)}t^{\gamma-1}{\rm d}t=\frac{\left(\frac2s(\gamma-1)+1\right)^{\frac s2}}\gamma u(x)^\gamma\geq\frac1\gamma u(x)u_L(x)^{\gamma-1}.
	\end{equation*}
	On the other hand, if $u(x)>L$, one has
	\begin{equation*}
		\begin{split}
			\Psi(u)(x)&=\int_0^L(\psi'(t))^{\frac s2}{\rm d}t+\int_L^{u(x)}(\psi'(t))^{\frac s2}{\rm d}t=\frac{\left(\frac2s(\gamma-1)+1\right)^{\frac s2}}\gamma L^\gamma+L^{\gamma-1}(u(x)-L)\\
			&\geq\frac{L^\gamma}\gamma+\frac{L^{\gamma-1}}\gamma(u(x)-L)=\frac1\gamma u(x)u_L(x)^{\gamma-1}.
		\end{split}
	\end{equation*}
	Hence, taking $\varepsilon=\frac1{2C}$ in \eqref{eqn:re4} and using the Sobolev embedding $W^{s,\frac2s}(\R^2)\hookrightarrow L^q(\R^2)$ with $q>\frac{2p}s$, one obtains
	$$
	\frac{S_q}{\gamma^{\frac2s}}\|u u_L^{\gamma-1}\|_q^{\frac2s}\leq C\|u u_L^{\gamma-1}\|_{\frac{2p}s}^{\frac2s}\,,
	$$
	where $S_q$ is the best constant of the embedding inequality (see \cite{CD} for related results). Letting $L\to+\infty$, we immediately deduce
	\begin{equation}\label{eqn:re10}
		\|u\|_{q\gamma}^{\frac2s}\leq C\frac{\gamma^{\frac2s}}{S_q}\|u \|_{\frac{2p\gamma}s}^{\frac2s}\,.
	\end{equation}
	Finally, buying the line of \cite[Lemma 21]{Bisci2022}, we get $\|u\|_\infty\leq C$.
\end{proof}

\noindent Next, in order to prove a crucial decay estimate, let us recall a standard lemma from \cite[Lemma 5.1]{Chen18}.
\begin{lemma}\label{Lem:mean}
	For $G(t)=|t|^{\frac2s-2}t$, it is well-known that by the mean value theorem, we have
	$$
	G(t_2)-G(t_1)=G'(\xi)(t_2-t_1)
	$$
	for some $\xi\in[t_1,t_2]$. Then there exists a constant $c_0>0$ independent of $t_1$ and $t_2$ such that
	$$
	|\xi|\geq c_0\max\{|t_1|,|t_2|\}\,.
	$$
\end{lemma}

\begin{lemma}\label{decay_polynomial}\label{Lem:G-decay}
	Let $u\in W^{s,\frac2s}(\R^2)$ be a positive solution of the Choquard equation \eqref{Chs}. Then there exist $C,R>0$ for which
	\begin{equation*}
		u(x)\leq\frac C{1+|x|^{\frac{7s}{2(2-s)}}}\qquad\mbox{for all}\ \ |x|>R\,.
	\end{equation*}
\end{lemma}

\begin{proof}
We first claim that
\begin{equation}\label{logF_estimate}
	\int_{\R^2}\ln\frac1{|x-y|}F(u(y))\,\dd y\leq\ln\frac2{|x|}\int_{\left\{|y|\leq\frac{|x|}2\right\}}F(u(y))\,\dd y+C_0,
\end{equation}
where $C_0$ is a constant independent of $x$. Indeed, denoting by 
$$
\Omega:=\left\{(x,y)\in\R^2\ :\ |y|>\frac{|x|}2,\ |x-y|<1\right\},
$$ 
we have
\begin{equation*}
	\begin{split}
		\int_{\R^2}\ln\frac1{|x-y|}F(u(y))\,\dd y&\leq\int_{|y|\leq\frac{|x|}2}\ln\frac2{|x|}F(u(y))\,\dd y+\int_\Omega\ln\frac1{|x-y|}F(u(y))\,\dd y\\
		&\leq\int_{|y|\leq\frac{|x|}2}\ln\frac2{|x|}F(u(y))\,\dd y
+C\int_{|x-y|<1}\ln\frac1{|x-y|}\,\dd y\,,
	\end{split}
\end{equation*}
with $C$ depending on $\|u\|_\infty$, from which \eqref{logF_estimate} follows. Hence, there exists $R_1>0$ such that
$$
\int_{\R^2}\ln\frac1{|x-y|}F(u(y))\,\dd y<0\quad\mbox{for all}\ |x|\geq R_1\,.
$$
Hence, taking $\varphi\in W^{s,\frac2s}(\R^2)$ such that $\varphi\geq0$ and $\textrm{supp}\,\varphi\subset B_{R_1}(0)^c$, we deduce
\begin{equation}\lab{eqn:q2}
	\aligned
		\int_{\R^2}\int_{\R^2}&\frac{|u(x)-u(y)|^{\frac2s-2}(u(x)-u(y))(\varphi(x)-\varphi(y))}{|x-y|^4}\,{\rm d}x\,{\rm d}y+\int_{\R^2}|u|^{\frac2s-2}u\varphi\,{\rm d}x\\
		&=\frac1{2\pi}\int_{\R^2}\left(\int_{\R^2}\ln\frac1{|x-y|}F(u(y))\,{\rm d}y\right)f(u(x))\varphi(x)\,{\rm d}x\leq0.
	\endaligned
\end{equation}
Defining now $w(x):=\big(1+|x|^{\frac{7s}{2(2-s)}}\big)^{-1}$, by \cite[Lemma 7.1]{Del20} there exists $R_2>0$ such that
$$|(-\Delta)^s_\frac2s w(x)|\leq\frac c{|x|^4}\quad\text{for}\quad |x|\geq R_2$$
and hence
\begin{equation}\lab{eqn:q3}
	(-\Delta)^s_\frac2s w+|w|^{\frac2s-2}w\geq0,\quad \text{for}\quad|x|\geq\widetilde R_2
\end{equation}
for a suitable $\widetilde R_2\geq R_2$. Note that $(-\Delta)^s_\frac2s w$ is pointwisely well-defined since $w\in C^{1,1}_{\loc}(\R^2)\cap L_{s,\frac2s}(\R^2)$. Since $u\in L^\infty(\R^2)$ by Lemma \ref{Lemma_infty}, one may find $C_1>0$, depending on $\|u\|_\infty$, for which
$$
\psi(x):=u(x)-C_1 w(x)\leq0\quad\quad\text{for}\quad |x|=R_3:=\max\{R_1,\widetilde R_2\}.
$$
Defining $\widetilde\psi:=\psi\chi_{\R^2\setminus B_{R_3}(0)}$ and taking $\widetilde\psi^+=\max\{\widetilde\psi,0\}\in W_0^{s,\frac2s}(B_{R_3}(0)^c)$ as a test function in \eqref{eqn:q2}, and noting that by homogeneity the inequality \eqref{eqn:q3} holds also for $\bar{w}:=C_1w$, we deduce that
\begin{equation}\lab{eqn:q4}
	\int_{\R^2}\int_{\R^2}\frac{h(x,y)(\widetilde\psi^+(x)-\widetilde\psi^+(y))}{|x-y|^4}\,{\rm d}x\,{\rm d}y+\int_{\R^2}\big(|u|^{\frac2s-2}u-|\bar{w}|^{\frac2s-2}\bar{w}\big)\widetilde\psi^+\,{\rm d}x\leq0\,,
\end{equation}
where
$$
h(x,y):=|u(x)-u(y)|^{\frac2s-2}(u(x)-u(y))-|\bar w(x)-\bar w(y)|^{\frac2s-2}(\bar w(x)-\bar w(y)).
$$
Note that by Lemma \ref{Lem:mean} one has
\begin{multline*}
\int_{\R^2}\int_{\R^2}\frac{h(x,y)(\widetilde\psi^+(x)-\widetilde\psi^+(y))}{|x-y|^4}\,{\rm d}x\,{\rm d}y\\
\geq\int_{\R^2}\int_{\R^2}\frac{|\bar h(x,y)|^{\frac2s-2}(\widetilde\psi(x)-\widetilde\psi(y)))(\widetilde\psi^+(x)-\widetilde\psi^+(y))}{|x-y|^4}\,{\rm d}x\,{\rm d}y\geq0\,,
\end{multline*}
where $\bar h(x,y)$ lies in the segment $[(u(x)-u(y)),(\bar w(x)-\bar w(y))]$. By \eqref{eqn:q4} this implies
$$
0\leq\int_{B_{R_3}(0)^c}\big(|u|^{\frac2s-2}u-|\bar w|^{\frac2s-2}\bar w\big)(u-\bar w)^+{\rm d}x\leq0\,.
$$
As a consequence we get 
$$
u(x)\leq C_1w(x)\leq \frac{C_1}{1+|x|^{\frac{7s}{2(2-s)}}}\quad\text{for}\quad |x|\geq R_3\,,
$$
as desired.
\end{proof}

\begin{lemma}\label{logFu}
	Let $u\in W^{s,\frac2s}(\R^2)$ be a positive solution of \eqref{Chs}. Then
	\begin{equation}\label{logFu_eq}
		\int_{\R^2}\ln(\e+|x|)F(u)\,{\rm d}x<+\infty\,.
	\end{equation}
\end{lemma}
\begin{proof}
	Lemma \ref{decay_polynomial} and ($f_2$) imply that there exists $R>0$ for which
	$$F(u(x))\leq C\left[1+|x|^{\frac{7s}{2(2-s)}}\right]^{-\frac2s}$$
	as $|x|>R$. Therefore, since $\|u\|_\infty\leq C$ by Lemma \ref{Lemma_infty}, one has
	\begin{equation*}
		\begin{split}
			\int_{\R^2}\ln(\e+|x|)F(u)\,{\rm d}x&=\int_{|x|<R}\ln(\e+|x|)F(u)\,{\rm d}x+\int_{|x|\geq R}\ln(\e+|x|)F(u)\,{\rm d}x\\
			&\leq C\int_{|x|<R}\ln(\e+|x|)\,{\rm d}x+C\int_{|x|\geq R}\ln(\e+|x|)\left[1+|x|^{\frac{7s}{2(2-s)}}\right]^{-\frac2s}\,{\rm d}x\\
			&\leq C+C\int_{|x|\geq R}|x|\cdot|x|^{-\frac7{2-s}}\,{\rm d}x<+\infty\,.
		\end{split}
	\end{equation*}
\end{proof}

\noindent We are now in the position to prove Theorem \ref{thm_equiv}.
\begin{proof}[Proof of Theorem \ref{thm_equiv}]
	We already proved that $u\in L^\infty(\R^2)$ in Lemma \ref{Lemma_infty}. To show that $\phi_u\in L_\frac12(\R^2)$ we compute as follows:
	\begin{equation*}
		\begin{split}
			2\pi\int_{\R^2}\frac{|\phi_u(x)|}{1+|x|^3}\,{\rm d}x&\leq\int_{\R^2}F(u(y))\left(\int_{\R^2}\left|\ln\left(\frac1{|x-y|}\right)\right|\frac1{1+|x|^3}\,{\rm d}x\right){\rm d}y\\
			&\leq\int_{\R^2}F(u(y))\left(\int_{|x-y|>1}\frac{\ln|x-y|}{1+|x|^3}\,{\rm d}x+\int_{|x-y|\leq1}\ln\frac1{|x-y|}\,{\rm d}x\right){\rm d}y\\
			&\leq\int_{\R^2}F(u(y))\,\dd y\int_{\R^2}\frac{\ln(1+|x|)}{1+|x|^3}\,\dd x+\int_{\R^2}\ln(1+|y|)F(u(y))\,\dd y\int_{\R^2}\frac{\dd x}{1+|x|^3}\\
			&\quad+\|F(u)\|_{L^1(\R^2)}\|\ln(\cdot)\|_{L^1(B_1(0))}<+\infty\,,
		\end{split}
	\end{equation*}
	using the elementary estimate $\ln|x-y|\leq\ln(1+|x|)+\ln(1+|y|)$ for $|x-y|>1$, and Lemma \ref{logFu}.
	
	\noindent In order to show that $(u,\phi_u)$ solves \eqref{SPs}, we follow the approach of \cite[Theorem 2.1]{BCT}. By \cite[Lemma 2.3]{H} the function
	$$\tv_u(x):=\frac1{2\pi}\int_{\R^2}\ln\left(\frac{1+|y|}{|x-y|}\right)F(u(y))\,{\rm d}y$$
	belongs to $H^1_{\loc}(\R^2)$ and solves $-\Delta\tv_u=F(u)$ in $\R^2$ in the sense of Definition \ref{sol_Poisson}. Moreover, by \cite[Lemma 2.4]{H} we know that all such solutions of $-\Delta\phi=F(u)$ in $\R^2$ are of the form $\phi=\tv_u+p$ with $p$ polynomial of degree at most one. We aim to prove that $\phi_u$ solves the Poisson equation in \eqref{SPs} by showing that $\tv_u-\phi_u$ is constant. Indeed,
	\begin{equation*}
		\begin{split}
			\tv_u(x)-\phi_u(x)&=\frac1{2\pi}\int_{\R^2}\left(\ln\left(\frac{1+|y|}{|x-y|}\right)-\ln\left(\frac1{|x-y|}\right)\right)F(u(y))\,{\rm d}y\\
			&=\frac1{2\pi}\int_{\R^2}\ln(1+|y|)F(u(y))\,{\rm d}y<+\infty
		\end{split}
	\end{equation*}
	by Lemma \ref{logFu}. Finally, we prove the logarithmic behaviour \eqref{asympt_phi} of $\phi_u$ at $\infty$ in the spirit of \cite[Proposition 2.3 (ii)]{CW}. One has
	\begin{equation*}
		2\pi\phi_u(x)+\ln|x|\int_{\R^2}F(u(y))\,\dd y=\int_{\R^2}\ln\frac{|x|}{|x-y|}F(u(y))\,\dd y=:-\int_{\R^2}l(x,y)F(u(y))\,\dd y\,,
	\end{equation*}
	where $l(x,y):=\ln\frac{|x-y|}{|x|}$. Note that for $|x|\to+\infty$ one has $l(x,y)\to0$ for any fixed $y\in\R^2$. Let now $|x|>1$, then if $|x-y|\geq\frac12$ one has
	$$\ln\frac12\leq\ln\frac12-\ln\frac1{|x|}\leq l(x,y)\chi_{|x-y|\geq\frac12}(y)\leq\ln\left(\frac{|x|+|x||y|}{|x|}\right)\leq\ln(1+|y|)\,.$$
	This yields
	\begin{equation*}
		\begin{split}
			\left|-\int_{|x-y|\geq\frac12}l(x,y)F(u(y))\,\dd y\right|&\leq\int_{\R^2}\max\left\{\ln2,\ln(1+|y|)\right\}F(u(y))\,\dd y\\
			&\leq\|F(u)\|_{L^1(\R^2)}\ln2+\int_{\R^2}\ln(1+|y|)F(u(y))\,\dd y<+\infty\,,
		\end{split}
	\end{equation*}
	by Lemma \ref{logFu}. By the dominated convergence theorem, we infer
	\begin{equation}\label{log_decay_1}
		\int_{|x-y|\geq\frac12}l(x,y)F(u(y))\,\dd y\to0\qquad\mbox{as}\ \,|x|\to+\infty\,.
	\end{equation}
	On the other hand, recalling that $|x|>1$, in the case $|x-y|\leq\frac12$, we have $|y|\geq\frac{|x|}2$, therefore
	\begin{equation*}
		0\leq\int_{|x-y|\leq\frac12}\ln|x|F(u(y))\,\dd y\leq\int_{\R^2\setminus B_{\frac{|x|}2}(0)}\ln(2(1+|y|))F(u(y))\,\dd y\to0
	\end{equation*}
	as $|x|\to+\infty$ for the same reason. Moreover, by ($f_2$) and Theorem \ref{Thm:m-theo} there exists $\beta>1$ for which $F(u)\in L^\beta(\R^2)$, and we have  that
	\begin{equation*}
		\begin{split}
			0\leq -\int_{|x-y|\leq\frac12}\ln|x-y|F(u(y))\,\dd y&=\int_{B_{\frac12}(0)}-\ln|z|F(u(x-z))\,\dd z\\
			&\leq\|\ln(\cdot)\|_{L^{\beta'}(B_\frac12(0))}\|F(u)\|_{L^\beta(B_\frac12(x))}\to0
		\end{split}
	\end{equation*}
	as $|x|\to+\infty$. Eventually we conclude that also
	$$\int_{|x-y|\leq\frac12}l(x,y)F(u(y))\,\dd y\to 0$$
	as $|x|\to+\infty$, which together with \eqref{log_decay_1} implies \eqref{asympt_phi}.
\end{proof}

\noindent Finally, we prove that solutions of \eqref{Chs} are H\"older continuous. Following \cite{dCKP,Iannizzotto16}, for all measurable $u:\R^2\to\R$, we define its $\big(s,\frac2s\big)-$\textit{non-local tail} centered at $x\in\R^2$ with radius $R>0$ as
$$\Tail(u;x,R):=R^2\left(\int_{\R^2\setminus B_R(x)}\frac{|u(y)|^{\frac2s-1}}{|x-y|^4}\,{\rm d}y\right)^{\frac s{2-s}}.$$

\begin{lemma}\label{Lemma_alpha,0}
	Let $u\in W^{s,\frac2s}(\R^2)$ be a positive solution of \eqref{Chs}, then $u\in C^\nu_{\loc}(\R^2)$ for some $\nu\in(0,1)$.
\end{lemma}
\begin{proof}
	First, as in the proof of Lemma \ref{Lemma_infty}, one has $(-\Delta)^s_{\frac2s}u\leq\bar C$ weakly in $\R^2$ for some constant $\bar C>0$. To prove the weak bound from below, take a test function $\varphi\geq0$ and estimate
	\begin{equation}\label{eqn:re1_holder}
		\begin{split}
			\int_{\R^2}\int_{\R^2}&\frac{|u(x)-u(y)|^{\frac2s-2}(u(x)-u(y))(\varphi(x)-\varphi(y))}{|x-y|^4}\,{\rm d}x\,{\rm d}y+\int_{\R^2}|u|^{\frac2s-2}u\varphi\,{\rm d}x\\
			&=\frac1{2\pi}\int_{\R^2}\left(\int_{\R^2}\ln\frac1{|x-y|}F(u(y))\,{\rm d}y\right)f(u(x))\varphi(x)\,{\rm d}x\\
			&\geq-\frac1{2\pi}\int_{|x-y|\geq1}\int_{\R^2}\ln|x-y|F(u(y))f(u(x))\varphi(x)\,{\rm d}x\,{\rm d}y\\
			&\geq-\frac1{2\pi}\int_{\R^2}\int_{\R^2}\left(\ln(1+|x|)+\ln(1+|y|)\right)F(u(y))f(u(x))\varphi(x)\,{\rm d}x\,{\rm d}y.
		\end{split}
	\end{equation}
	Since $u\in L^\infty(\R^2)$ and by the decay estimate of Lemma \ref{Lem:G-decay} we have
	\begin{equation*}
		\begin{split}
			\int_{\R^2}\int_{\R^2}&\ln(1+|x|)F(u(y))f(u(x))\varphi(x)\,{\rm d}x\,{\rm d}y\\
			&\leq C\int_{\R^2}\frac{\ln(1+|x|)\varphi(x)}{\left(1+|x|^\frac{7s}{2(2-s)}\right)^{\frac2s-1}}\,\dd x\int_{\R^2}\frac{\dd y}{\left(1+|y|^\frac{7s}{2(2-s)}\right)^\frac2s}\\
			&\leq C\int_{\R^2}\frac{\ln(1+|x|)\varphi(x)}{1+|x|^\frac72}\,\dd x\int_{\R^2}\frac{\dd y}{1+|y|^\frac7{2-s}}\leq C\int_{\R^2}\varphi\,,
		\end{split}
	\end{equation*}
	and analogously
	\begin{equation*}
		\begin{split}
			\int_{\R^2}\int_{\R^2}&\ln(1+|y|)F(u(y))f(u(x))\varphi(x)\,{\rm d}x\,{\rm d}y\leq C\int_{\R^2}\frac{\varphi(x)}{1+|x|^\frac72}\,\dd x\int_{\R^2}\frac{\ln(1+|y|)}{1+|y|^\frac7{2-s}}\,\dd x\leq C\int_{\R^2}\varphi.
		\end{split}
	\end{equation*}
	Since $\int_{\R^2}|u|^{\frac2s-2}u\varphi\leq C\int_{\R^2}\varphi$, this yields $(-\Delta)^s_{\frac2s}u\geq -C$ weakly. Finally, $\Tail(u;x,R)$ is bounded uniformly in $(x,R)$, since
	\begin{equation}\label{Tail}
		\Tail(u;x,R)\leq R^2\|u\|_\infty\left(\int_{\R^2\setminus B_R(0)}\frac{{\rm d}y}{|y|^4}\right)^{\frac s{2-s}}\leq C\,.
	\end{equation}
	Then, by \cite[Corollary 5.5]{Iannizzotto16} there exists a universal constant $\tC$ and $\nu\in(0,1)$ such that for all $x_0\in\R^2$ one has
	\begin{equation*}
		[u]_{C^\nu(B_R(x_0))}\leq\tC R^{-\nu}\left((\bar CR^2)^{\frac s{2-s}}+\|u\|_{L^\infty(B_{2R}(x_0))}+\Tail(u;x_0,2R)\right)\leq C(R)
	\end{equation*}
	with $C(R)$ independent of the choice of $x_0$. This readily implies $u\in C^\nu_{\loc}(\R^2)$.
\end{proof}

\section{Symmetry of positive solutions: proof of Theorem \ref{Thm:symmetric}}\label{Sec_MVP}

\noindent As a consequence of Theorem \ref{thm_equiv}, we have that solutions $u$ of \eqref{Chs} (which are positive, see Remark \ref{Rmk_positivity}) correspond to solutions $(u,\phi)$ of \eqref{SPs} which enjoy the following 
\begin{equation}\label{eqn:q1+0}
	u\in L^\infty(\R^2)\quad\quad\text{and}\quad\quad \phi(x)=\phi_u(x)=\left(I_2\ast F(u)\right)(x)\to-\infty\quad\text{as}\ \,|x|\to+\infty.
\end{equation}
\noindent Next we are going to prove Theorem \ref{Thm:symmetric} and for this purpose we rely on the method of moving planes, which has been adapted to the $p$-fractional context by Chen and Li in \cite{Chen18}, provided the operator $(-\Delta)^s_\frac2s$ is pointwisley defined. For this reason we require $u\in C_{\loc}^{1,1}(\R^2)\cap W^{s,\frac2s}(\R^2)$, see \cite[Lemma 5.2]{Chen18} and also \cite{WC}.

\medskip

\noindent For $\beta\in\R$, let us set
$$
\Sigma_\beta:=\{x\in\R^2:x_1<\beta\}\ \quad\mbox{and}\ \quad\partial \Sigma_\beta=\{x\in\R^2:x_1=\beta\}.
$$
Moreover, for any $x\in\R^2$, denote by $x^\beta$ the reflection of $x$ with respect to $\partial \Sigma_\beta$, that is $x^\beta=(2\beta-x_1,x_2)$. Set also
$$
u^{\beta}(x):=u(x^{\beta})\ \quad\mbox{and}\ \quad\phi^{\beta}(x):=\phi(x^{\beta})\quad\ \text{for}\ x\in\R^2.
$$
Note that $\phi$ is a continuous function on $\R^2$ by ($f_1$) and Lemma \ref{Lem:G-decay}. Defining
$$u_\beta:=u^{\beta}-u\quad\ \mbox{and}\ \quad\phi_\beta:=\phi^{\beta}-\phi\,,$$
the following holds
\begin{equation}\label{eqn:mod-2}
	(-\Delta)^s_{\frac2s} u^\beta-(-\Delta)^s_{\frac2s} u +|u^\beta|^{\frac2s-2}u^\beta-|u|^{\frac2s-2}u=\phi_\beta f(u^\beta)+\phi h_\beta u_\beta,\quad x\in \Sigma_\beta,
\end{equation}
where
$$
h_\beta(x):=\left\{
  \begin{array}{ll}
    \frac{f(u^\beta(x))-f(u(x))}{u_\beta(x)}, &  u_\beta(x)\not=0,\\
    f'(u(x)), & u_\beta(x)=0.
  \end{array}
\right.
$$
Recalling that $f\in C^1(\R)$, there exists $C=C(u)>0$ such that $\|h_\beta\|_{L^\infty(\Sigma_\beta)}\leq C$ for any $\beta\in\R$. Moreover,
\begin{equation}\label{eqn:diff}
	-\Delta\phi_\beta=K_\beta(x)u_\beta,\quad x\in \Sigma_\beta,
\end{equation}
where
$$
K_\beta(x):=\left\{
  \begin{array}{ll}
    \frac{F(u^\beta(x))-F(u(x))}{u_\beta(x)}, &  u_\beta(x)\not=0,\\
    f(u(x)), & u_\beta(x)=0.
  \end{array}
\right.
$$
From the definition of $\phi$, we deduce that
\begin{equation}\lab{eqn:phi-expression}
	\phi_\beta(x)=\int_{\Sigma_\beta}\ln\frac{|x-y^\beta|}{|x-y|}K_{\beta}(y)u_\beta(y)\,{\rm d}y\,.
\end{equation}
Since $\frac{|x-y^\beta|}{|x-y|}>1$ for any $x,y\in \Sigma_\beta$, we have $\phi_\beta\geq 0$ in $\Sigma_\beta$ if $u_\beta\geq 0$ in $\Sigma_\beta$ for every $\beta\in\R$. In what follows, for a function $v$ we denote by $v^-:=\min\{v,0\}$ its negative part.
Notice that $\phi^-$ is a nonpositive function with this convention. Actually we have the following lemma which is inspired by \cite[Theorem 2.1]{Chen18}.

\begin{lemma}\label{Lem:w-guji3}
	If $\beta\in\R$ is such that $u_\beta\geq0$ in $\Sigma_\beta$, then also $\phi_\beta\geq0$ on $\Sigma_\beta$. Furthermore,
	either $u_\beta\equiv0\equiv \phi_\beta$ or $u_\beta>0, \phi_\beta>0$ on $\Sigma_\beta$.
\end{lemma}
\begin{proof}
	From the expression of $\phi_\beta$ in \eqref{eqn:phi-expression} we know that $u_\beta\geq 0$ in $\Sigma_\beta$ implies $\phi_\beta\geq0$ in $\Sigma_\beta$. In particular, $\phi_\beta>0$ in $\Sigma_\beta$ if $u_\beta\not\equiv0$. Conversely, if $\phi_\beta\not\equiv0$, then also $u_\beta\not\equiv0$ again using \eqref{eqn:phi-expression}. We \textit{claim}:	\begin{equation}\label{eqn:positive}
	u_\beta>0\quad\quad\text{for}\quad x\in \Sigma_\beta\,.
	\end{equation}
	Indeed, if not, there exists $x_0\in \Sigma_\beta$ such that
	$u_\beta(x_0)=\min_{\Sigma_\beta}u_\beta=0$. Let $G(t)=|t|^{\frac2s-2}t$, then $G(t)$ is a strictly increasing function, and $G'(t)=(\frac2s-1)t^{\frac2s-2}\geq0$.
	On the one hand, a direct computation using the definition of the $(s,\frac2s)$-fractional operator yields
	\begin{equation}\label{eqn:guji}
		\begin{split}
			(-&\Delta)^s_{\frac2s}u^\beta(x_0)-(-\Delta)^s_{\frac2s}u(x_0)+|u^\beta(x_0)|^{\frac2s-2}u^\beta(x_0)-|u(x_0)|^{\frac2s-2}u(x_0)\\
			&=C\cdot PV\int_{\R^2}\frac{G(u^\beta(x_0)-u^\beta(y))-G(u(x_0)-u(y))}{|x-y|^4}\,{\rm d}y\\
			&=C\cdot PV\int_{\Sigma_\beta}\frac{G(u^\beta(x_0)-u^\beta(y))-G(u(x_0)-u(y))}{|x-y|^4}\,{\rm d}y\\
			&\quad+C\cdot PV\int_{\Sigma_\beta}\frac{G(u^\beta(x_0)-u(y))-G(u(x_0)-u(y^\beta))}{|x-y^\beta|^4}\,{\rm d}y\\
			&=C\cdot PV\int_{\Sigma_\beta}\bigg[\frac{1}{|x-y|^4}-\frac1{|x-y^\beta|^4}\bigg]\left(G(u^\beta(x_0)-u^\beta(y))-G(u(x_0)-u(y))\right){\rm d}y\\
			&\quad+C\cdot PV\int_{\Sigma_\beta}\frac{G'(\xi(y))+G'(\eta(y))}{|x-y|^4}\,u_\beta(x_0)\,{\rm d}y\\
			&=-C\cdot PV\int_{\Sigma_\beta}\bigg[\frac1{|x-y|^4}-\frac{1}{|x-y^\beta|^4}\bigg]G'(\theta(y))u_\beta(y)\,{\rm d}y\leq0\,,
		\end{split}
	\end{equation}
	where $\xi(y)$ lies between $u^\beta(x_0)-u^\beta(y)$ and $u(x_0)-u^\beta(y)$, $\eta(y)$ between $u^\beta(x_0)-u(y)$ and $u(x_0)-u(y)$, and $\theta(y)$ between $u^\beta(x_0)-u^\beta(y)$ and $u(x_0)-u(y)$. Here we have used the mean value theorem in the third and fourth identities, and the fact that
	$$
	\frac{1}{|x-y|^4}>\frac{1}{|x-y^\beta|^4},\quad\forall x,y\in \Sigma_\beta.
	$$
	On the other hand, by assumption $(f_1)$ and since $w_\beta(x_0)>0$, observe that
	$$
	\phi_\beta(x_0) f(u^\beta(x_0))+\phi(x_0)h_\beta(x_0) u_\beta(x_0)=\phi_\beta(x_0) f(u^\beta(x_0))>0\,,
	$$
	which implies by \eqref{eqn:mod-2} that
	\begin{equation}\label{eqn:big0}
	(-\Delta)^s_{\frac2s} u^\beta(x_0)-(-\Delta)^s_{\frac2s} u (x_0) +|u^\beta(x_0)|^{\frac2s-2}u^\beta-|u(x_0)|^{\frac2s-2}u(x_0)>0\,,
	\end{equation}
	which contradicts \eqref{eqn:guji}. 
\end{proof}

\noindent Let us next recall a maximum principle for the $p-$fractional laplacian and a key boundary estimate lemma for anti-symmetric functions, both established in \cite{Chen18}. These results play an important role in carrying out the method of moving planes.

\begin{lemma}\label{Lem:MP-Anti}
	(A maximum principle for anti-symmetric functions). Let $\Omega$ be a bounded domain in $\Sigma_\beta$. Assume that $u\in C_{loc}^{1,1}(\R^2)\cap W^{s,\frac2s}(\R^2)$. If
	\begin{equation}\lab{eqn:q1+}
		\begin{cases}
			(-\Delta)^s_{\frac2s} u^\beta(x)\geq(-\Delta)^s_{\frac2s}u(x)\ \  &x\in\Omega, \\
			u^\beta(x)\geq u(x) & x\in\Sigma_{\beta}\setminus{\Omega},
		\end{cases}
	\end{equation}
	then $u_\beta(x):=u^\beta(x)-u(x)\geq0$ in $\Omega$. If $u_\beta=0$ at some point in $\Omega$, then	$u_\beta(x)=0$ a.e. $x\in\R^2$. These conclusions hold for an unbounded region $\Omega$ if we further assume that
	$$\lim_{|x|\to+\infty}u_\beta(x)\geq0\,.$$
\end{lemma}

\begin{lemma}\label{Lem:bound}
	(A key boundary estimate). Assume that $u_{\beta_0}>0$ for $x\in \Sigma_{\beta_0}$.	Suppose $\beta_k\rightarrow\beta_0$, and $x^{k}\in \Sigma_{\beta_k}$, such that
	$$
	u_{\beta_k}(x^k)=\min_{\Sigma_{\beta_k}}u_{\beta_k}\leq0\quad\quad\text{and}\quad\quad x^k\rightarrow x_0\in\partial\Sigma_{\beta_0}\,.
	$$
	Let $\delta_k:={\rm dist}(x^k,\partial\Sigma_{\beta_k})\equiv|\beta_k-x_1^k|$. Then
	$$
	\limsup_{\delta_k\to0}\frac1{\delta_k}\left\{(-\Delta)^s_{\frac2s}u^\beta(x^k)-(-\Delta)^s_{\frac2s}u(x^k)\right\}<0\,.
	$$
\end{lemma}

\begin{lemma}\label{Lem:start}
	There exists $\beta<0$ such that
	\begin{equation}\label{eqn:posi1}
		u_\beta(x)\geq0,\quad\ \mbox{for all}\ \ x\in\Sigma_\beta\,.
	\end{equation}
\end{lemma}
\begin{proof}
	Note that $u$ satisfies \eqref{eqn:mod-2}. Suppose on the contrary that \eqref{eqn:posi1} is violated. Then, there exists $\beta_0<0$ such that for all $\beta<\beta_0$ one is always able to find $x_*\in\Sigma_\beta$ for which
	$$
	u_\beta(x_*)=\min\limits_{\Sigma_\beta}u_\beta<0
	$$
	holds First, observe that
	\begin{equation}\label{eqn:posi2+}
		0\leq \ln\frac{|x_*-y^\beta|}{|x_*-y|}\leq \ln\big(1+\frac{|y-y^\beta|}{|x_*-y|}\big)\leq\frac{|y-y^\beta|}{|x_*-y|}=\frac{2(\beta-y_1)}{|x_*-y|}
	\end{equation}
	for any $x,y\in \Sigma_\beta$. Set
	$$
	M_\beta:=\{u\in\Sigma_\beta\,:\, u_\beta<0 \}\,.
	$$
	Recalling that, by assumption ($f_1$), $K_\beta$ is non-negative function. Based on the above facts, we deduce from the expression of $\phi_\beta$ in \eqref{eqn:phi-expression} that
	\begin{equation}\label{eqn:posi3}
		\begin{split}
			\phi_\beta(x_*):&=\int_{\Sigma_\beta}\ln\frac{|x_*-y^\beta|}{|x_*-y|}K_{\beta}(y)u_\beta(y)\,{\rm d}y\\
			&\geq u_\beta(x_*)\int_{M_\beta}\frac{2(\beta-y_1)}{|x_*-y|}K_{\beta}(y)\,{\rm d}y\\
			&\geq u_\beta(x_*)\int_{M_\beta}\frac{2(\beta-y_1)}{|x_*-y|}f(u(y))\,{\rm d}y\\
			&\geq u_\beta(x_*)\int_{M_\beta\cap\{|x_*-y|\}>1}\frac{2|y|}{|x_*-y|}f(u(y))\,{\rm d}y\\
			&\quad+u_\beta(x_*)\int_{M_\beta\cap\{|x_*-y|\}\leq1}\frac{2(\beta-y_1)}{|x_*-y|}f(u(y))\,{\rm d}y\\
			&\geq u_\beta(x_*)\int_{\Sigma_\beta}2|y|f(u(y)){\rm d}y
			\\
			&\quad+u_\beta(x_*)\bigg(\int_{\{|x_*-y|\leq1\}}\frac{{\rm d}y}{|x_*-y|^\frac32}\bigg)^{\frac23}\bigg(\int_{\Sigma_\beta}|2(\beta-y_1)f(u(y))|^3\,{\rm d}y\bigg)^{\frac13}\\
			&=:\,u_\beta(x_*)c_\beta\,.
		\end{split}
	\end{equation}
	It is easy to check that $c_\beta$ is finite and $c_\beta\rightarrow0$ as $\beta\rightarrow-\infty$. Indeed $|y|f(u)\in L^1(\R^2)$ by ($f_2$) by Theorem \ref{Thm:m-theo}, and the second integral is finite. Moreover, for $|\beta|$ large enough in $B_{|\beta|}(0)^c$ one has $|u(x)|\leq C\left(1+|x|^{\frac{7s}{2(2-s)}}\right)$ by Lemma \ref{decay_polynomial}, and using ($f_1$) we obtain 
	\begin{equation*}
		\begin{split}
			\int_{\Sigma_\beta}|2(\beta-y_1)f(u(y))|^3{\rm d}y\leq C\int_{B_{|\beta|}(0)^c}|\beta-y_1|^3|y|^{\frac{21}2}{\rm d}y\leq C\beta^{-\frac{11}2}\to0\quad\mbox{as}\ \beta\to-\infty\,.
		\end{split}
	\end{equation*}  Moreover, by assumption ($f_1$), for any fixed $\varepsilon>0$ and $\beta$ sufficiently negative,
	\begin{equation}\label{eqn:posi4}
		f(u^\beta(x_*))\leq f(u(x_*))\leq \varepsilon|u(x_*)|^{\frac2s-2}u(x_*)\,.
	\end{equation}
	Thanks to \eqref{eqn:q1+0} and the fact that $h_\beta$ is a non-negative function, for sufficiently negative $\beta$,
	\begin{equation}\label{eqn:posi4+}
		\phi(x_*) h_\beta(x_*) u_\beta(x_*)\geq0\,.
	\end{equation}
	Combining \eqref{eqn:posi3}, \eqref{eqn:posi4}, and \eqref{eqn:posi4+} yields
	\begin{equation}\label{eqn:posi5}
		\phi_\beta(x_*)f(u^\beta(x_*))+\phi(x_*) h_\beta(x_*) u_\beta(x_*)\geq c_\beta |u(x_*)|^{\frac2s-2}u(x_*)u_\beta(x_*)\,.
	\end{equation}
It follows from Lemma \ref{Lem:mean} that there exists $c_o>0$ such that
	\begin{equation}\label{eqn:posi6}
		|u^\beta(x_*)|^{\frac2s-2}u^\beta(x_*)-|u(x_*)|^{\frac2s-2}u(x_*)\geq c_o|u(x_*)|^{\frac2s-2}u_\beta(x_*)\,.
	\end{equation}By using \eqref{eqn:posi5} and \eqref{eqn:posi6} in \eqref{eqn:mod-2} and $\beta$ sufficiently negative so that $c_\beta u(x_*)<c_o$, we get
	\begin{equation}\label{eqn:posi7}
		(-\Delta)^s_{\frac2s} u^\beta (x_*)-(-\Delta)^s_{\frac2s} u(x_*)>0\,.
	\end{equation}
Next, similar computations as in \eqref{eqn:guji} yield
	\begin{equation}\label{eqn:posi8}
		\aligned
		&(-\Delta)^s_{\frac2s} u^\beta(x_*)-(-\Delta)^s_{\frac2s} u(x_*)\\
		&=\,C\cdot PV\int_{\Sigma_\beta}\bigg[\frac{1}{|x-y|^4}-\frac{1}{|x-y^\beta|^4}\bigg]\left(G(u^\beta(x_*)-u^\beta(y))-G(u(x_*)-u(y))\right){\rm d}y\\
		&\quad+C\cdot PV\int_{\Sigma_\beta}\frac{G'(\xi(y))+G'(\eta(y))}{|x-y|^4}\,u_\beta(x_*)\,{\rm d}y\\
		&\leq\,C\cdot PV\int_{\Sigma_\beta}\bigg[\frac{1}{|x-y|^4}-\frac{1}{|x-y^\beta|^4}\bigg]\left(G(u^\beta(x_*)-u^\beta(y))-G(u(x_*)-u(y))\right){\rm d}y\,.\\
		\endaligned
	\end{equation}
	Since $G(t)=|t|^{\frac2s-2}t$ is a strictly increasing function, and
	$$
	u^\beta(x_*)-u^\beta(y)-u(x_*)+u(y)=u_\beta(x_*)-u_\beta(y)\leq0\,,
	$$
	we infer that
	$$
	G(u^\beta(x_*)-u^\beta(y))-G(u(x_*)-u(y))\leq0\,.
	$$
	This implies by \eqref{eqn:posi8} that
	$$
	(-\Delta)^s_{\frac2s} u^\beta(x_0)-(-\Delta)^s_{\frac2s} u(x_*)\leq0\,,
	$$
	which contradicts \eqref{eqn:posi7}. 	
\end{proof}

\begin{proof}[Proof of Theorem \ref{Thm:symmetric}]
	So far, Lemma \ref{Lem:start} provides a starting point to move the plane $\partial\Sigma_\beta$. Now, let us move the plane to the right, as long as \eqref{eqn:posi1} holds, up to some limiting position. In particular, define
	\begin{equation}\label{beta_0}
		\beta_0:=\sup\{\beta\in\R\,|\,u_\mu(x)\geq0,\,x\in\Sigma_\mu,\,\mu\leq\beta\}\,.
	\end{equation}

	\noindent Let us divide the proof into two steps:
	
	\noindent \emph{Step 1.} Let us prove that $u$ is symmetric about the limiting plane $\Sigma_{\beta_0}$, that is,
	\begin{equation}\label{eqn:posi9}
		\aligned
		u_{\beta_0}\equiv0,\quad\ \mbox{for all}\ \ x\in\Sigma_{\beta_0}\,.
		\endaligned
	\end{equation}
	By contradiction, assume that \eqref{eqn:posi9} is false, that is to say, $u_{\beta_0}(x)\geq0$ and $u_{\beta_0}(x)\neq0$ for some $x\in\Sigma_{\beta_0}$. Moreover, it follows from \eqref{eqn:mod-2} that
	\begin{equation}\label{eqn:posi10}
		(-\Delta)^s_{\frac2s} u^{\beta_0}-(-\Delta)^s_{\frac2s} u +|u^{\beta_0}|^{\frac2s-2}u^{\beta_0}-|u|^{\frac2s-2}u-\phi h_{\beta_0} u_{\beta_0}\geq0,\quad x\in \Sigma_{\beta_0}\,.
	\end{equation}
	So, by the strong maximum principle (see Lemma \ref{Lem:MP-Anti}) together with the fact that the map $t\mapsto|t|^{\frac2s-2}t$ is increasing and $t\mapsto\phi h_{\beta_0}t$ is linear, we have $u_{\beta_0}>0$ and hence by \eqref{eqn:phi-expression} also
	$$
	\phi_{\beta_0}(x)>0, \quad \forall x\in \Sigma_{\beta_0}
	$$
	by Lemma \ref{Lem:w-guji3}. Now, according to the definition of $\beta_0$, there exists a sequence $\beta_n\searrow \beta_0$, and $x^n\in\Sigma_{\beta_n}$ such that
	\begin{equation}\label{eqn:posi11}
		u_{\beta_n}(x^n)=\min\limits_{\Sigma_{\beta_n}}u_{\beta_n}<0,\quad\text{and}\quad \nabla u_{\beta_n}(x^n)=0\,.
	\end{equation}
	Remind that $u\in C^1(\R^2)$. Up to subsequence, we \textit{claim}:
	\begin{equation}\label{eqn:posi11+}
		x^n\rightarrow x^*\quad\text{as}\quad n\rightarrow+\infty\,.
	\end{equation}
	Indeed, let $B_R:=B_R(0)$ for $R>0$, by \eqref{eqn:q1+0}, we can choose $R>1$ large enough such that
	\begin{equation}\label{eqn:posi12}
		\phi\leq0\quad\text{in}\,\, \Sigma_{\beta_n}\setminus{B_R}\quad \text{for\,any\,}\beta_n
	\end{equation}
	and
	\begin{equation}\label{eqn:posi12+0}
		\aligned
		&\bigg(\int_{|x^n-y|\leq1}\frac{{\rm d}y}{|x^n-y|^\frac32}\bigg)^\frac23\left(\int_{\Sigma_{\beta_n}\setminus{B_R}}|2({\beta_n}-y_1)f(u(y))|^3{\rm d}y\right)^\frac13<\frac{c_o}4\,,\\
		&\int_{\Sigma_{\beta_n}\setminus{B_R}}2|y|f(u(y))\,{\rm d}y<\frac{c_o}4\quad\text{in}\,\, \Sigma_{\beta_n}\setminus{B_R}\quad \text{for\,any\,}\beta_n\,,
		\endaligned
	\end{equation}
	where $c_o$ has appeared in \eqref{eqn:posi6}. Assume by contradiction that there exists $N>0$ such that $x^n\in \Sigma_{\beta_n}\setminus{B_R}$ for all $n>N$. Then by \eqref{eqn:posi12} and the definition of $h_\beta$, one has for $n>N$
	\begin{equation}\label{eqn:posi12+}
		\phi(x^n) h_{\beta_n}(x^n)u_{\beta_n}(x^n)\geq0\,.
	\end{equation}
	Moreover, by assumption ($f_1$), for any fixed $\varepsilon>0$, taking $N$ large enough,
	\begin{equation}\label{eqn:posi13+0}
		f(u^{\beta_n}(x^n))\leq f(u(x^n))\leq \varepsilon|u(x^n)|^{\frac2s-2}u(x^n)
	\end{equation}
	holds for $n>N$. Note that for $y\in B_R$ one has
	$$
	\ln\frac{|x^n-y^{\beta_n}|}{|x^n-y|}\rightarrow 0\ \quad\text{as}\,\, n\rightarrow\infty\,.
	$$
	Set
	$$M_{\beta_n}:=\{u\in\Sigma_{\beta_n}\,:\, u_{\beta_n}<0 \}\,.$$
	We deduce from \eqref{eqn:phi-expression}, Lemma \ref{Lem:G-decay}, \eqref{eqn:posi2+}, and \eqref{eqn:posi12+0} that for $n>N$ with $N$ large enough
	\begin{equation}\label{eqn:posi13}
		\begin{split}
			\phi_{\beta_n}(x^n):&=\int_{\Sigma_{\beta_n}}\ln\frac{|x^n-y^{\beta_n}|}{|x^n-y|}K_{{\beta_n}}(y)u_{\beta_n}(y)\,{\rm d}y\\
			&\geq\int_{M_{\beta_n}\cap B_R} \ln\frac{|x^n-y^{\beta_n}|}{|x^n-y|}K_{{\beta_n}}(y)u_{\beta_n}(y)\,{\rm d}y\\
			&\quad+\int_{M_{\beta_n}\setminus{B_R}}\frac{2({\beta_n}-y_1)}{|x^n-y|}K_{{\beta_n}}(y)u_{\beta_n}(y)\,{\rm d}y\\
			&\geq u_{\beta_n}(x^n)\int_{M_{\beta_n}\setminus{B_R}}\frac{2({\beta_n}-y_1)}{|x^n-y|}f(u(y))\,{\rm d}y+o_n(1)\cdot u_{\beta_n}(x^n)\\
			&\geq u_{\beta_n}(x^n)\int_{M_{\beta_n}\cap\{|x^n-y|\geq1\}\setminus{B_R}}\frac{2|y|}{|x^n-y|}f(u(y))\,{\rm d}y\\
			&\quad+u_{\beta_n}(x^n)\int_{M_{\beta_n}\cap\{|x^n-y|\leq1\}\setminus{B_R}}\frac{2({\beta_n}-y_1)}{|x^n-y|}f(u(y))\,{\rm d}y+o_n(1)\cdot u_{\beta_n}(x^n)\\
			&\geq u_{\beta_n}(x^n)\bigg(\int_{|x^n-y|\leq1}\frac{{\rm d}y}{|x^n-y|^\frac32}\bigg)^\frac23\bigg(\int_{\Sigma_{\beta_n}\setminus{B_R}}|2({\beta_n}-y_1)f(u(y))|^3{\rm d}y\bigg)^\frac13\\
			&\quad+u_{\beta_n}(x^n)\int_{\Sigma_{\beta_n}\setminus{B_R}}2|y|f(u(y))\,{\rm d}y+o_n(1)\cdot u_{\beta_n}(x^n)\\
			&\geq\frac{c_o}2u_{\beta_n}(x^n)\,.
		\end{split}
	\end{equation}
	On the other hand, as for \eqref{eqn:posi6}, it follows from Lemma \ref{Lem:mean} that for $n>N$ with $N$ large enough
	\begin{equation}\label{eqn:posi14}
		|u(x^n)|^{\frac2s-2}u(x^n)-|u^{\beta_n}(x^n)|^{\frac2s-2}u^{\beta_n}(x^n)\geq c_o|u(x^n)|^{\frac2s-2}u_{\beta_n}(x^n)\,.
	\end{equation}
	Combining \eqref{eqn:posi12+}, \eqref{eqn:posi13+0} and \eqref{eqn:posi13}, we have
	$$
	\phi_{\beta_n}(x^n)f(u^{\beta_n}(x^n))+\phi(x^n) h_{\beta_n}(x^n) u_{\beta_n}(x^n)\geq \varepsilon \frac{c_o}{2}|u(x^n)|^{\frac2s-2}u(x^n)u_{\beta_n}(x^n)\,.
	$$
	Take $\varepsilon>0$ sufficiently small such that $\varepsilon u(x^n)<2$, when $n$ is large enough.
	Hence, from \eqref{eqn:posi14} and \eqref{eqn:mod-2} we immediately deduce that for $n>N$ with $N$ large enough
	$$
	(-\Delta)^s_{\frac2s} u_{\beta_n}(x^n)-(-\Delta)^s_{\frac2s} u(x^n)>0\,.
	$$
	Moreover, arguing as in Lemma \ref{Lem:start}, one may deduce
	$$
	(-\Delta)^s_{\frac2s} u_{\beta_n}(x^n)-(-\Delta)^s_{\frac2s} u(x^n)<0\,.
	$$
	This is a contradiction. Therefore, $(x^n)_n$ must be bounded. The claim \eqref{eqn:posi11+} holds true.  Thus, by \eqref{eqn:posi11} we have
	\begin{equation}\label{eqn:posi15}
		u_{\beta_0}(x^*)=\min\limits_{\Sigma_{\beta_0}}u_{\beta_0}\leq0,\quad\text{and}\quad \nabla u_{\beta_0}(x^*)=0\,.
	\end{equation}
	By the definitions of $\delta_n, h_{\beta_n}$ and recalling that $\phi$ is a continuous function, it follows that
	\begin{equation}\label{eqn:posi16}
		\aligned
		\liminf\limits_{n\rightarrow\infty}\frac{\phi(x^n)h_{\beta_n}(x^n)u_{\beta_n}(x^n)}{\delta_n}=\phi(x^*)h_{\beta_0}(x^*)\liminf\limits_{n\rightarrow\infty}\frac{u_\beta(x^n)}{\delta_n}=0\,,
		\endaligned
	\end{equation}
	where $\delta_n:=|\beta_n-x^n_1|$. By the definitions of $K_{\beta}$ and $M_\beta$, Lemma \ref{Lem:G-decay}, \eqref{eqn:posi2+}, and ($f_1$)-($f_2$), there exists $C>0$ such that for $R>1$
	\begin{equation}\label{eqn:posi17}
		\begin{split}
			\phi_{\beta_n}(&x^n)=\int_{\Sigma_{\beta_n}}\ln\frac{|x^n-y^{\beta_n}|}{|x^n-y|}K_{{\beta_n}}(y)u_{\beta_n}(y)\,{\rm d}y\\
			&\geq u_{\beta_n}(x^n)\bigg(\int_{M_{\beta_n}\cap B_R} \frac{2({\beta_n}-y_1)}{|x^n-y|}f(u(y))\,{\rm d}y+\int_{M_{\beta_n}\setminus{B_R}}\frac{2({\beta_n}-y_1)}{|x^n-y|}f(u(y))\,{\rm d}y\bigg)\\
			&\geq u_{\beta_n}(x^n)\bigg(\int_{M_{\beta_n}\cap B_R}\frac{2|y|}{|x^n-y|}f(u(y))\,{\rm d}y+\int_{M_{\beta_n}\setminus{B_R}}2|y|f(u(y))\,{\rm d}y\bigg)\\
			&\geq Cu_{\beta_n}(x^n)\,.
		\end{split}
	\end{equation}
	Hence
	$$
	\liminf\limits_{n\rightarrow\infty}\frac{\phi_{\beta_n}(x^n)f(u^{\beta_n}(x_n))}{\delta_n}\geq\liminf\limits_{n\rightarrow\infty}\frac{Cu_{\beta_n}(x^n)f(u^{\beta_n}(x_n))}{\delta_n}=0\,,
	$$
	which implies, together with \eqref{eqn:posi17} and \eqref{eqn:mod-2} for $\beta=\beta_n$ and $x=x^n$, that
	$$
	\liminf\limits_{\delta_n\to 0^+}\frac1{\delta_n}\left\{(-\Delta)^s_{\frac2s} u^\beta(x^n)-(-\Delta)^s_{\frac2s}u(x^n)\right\}\geq0\,.
	$$
	This contradicts Lemma \ref{Lem:bound} and therefore \eqref{eqn:posi9} holds true.

\medskip	

	\noindent \emph{Step 2.} We next complete the proof showing that $u$ is radially symmetric.  Recalling Lemma \ref{Lem:G-decay} and the definition of $\beta_0$ in \eqref{beta_0}, we first have $\beta_0<\infty$. It follows from Lemma \ref{Lem:bound} that $\beta_0>-\infty$. According to Lemma \ref{Lem:w-guji3} and Step 1, we get $u_{\beta_0}\equiv0$ and $\phi_{\beta_0}\equiv0$. By using the same argument for the second coordinate direction $x_2$, we can find $\beta_2\in\R$ such that $u_{\beta_2}\equiv0$ and $\phi_{\beta_2}\equiv0$. Consider $\beta=(\beta_0,\beta_2)$, then $\tilde u(x):=u(x-\beta)$ and $\tilde\phi(x):=\phi(x-\beta)$ is a solution of equation \eqref{eqn:mod-2}. By invariance under translation, we may assume that $\tilde{u}(x)=\tilde{u}(-x)$ and $\tilde{\phi}(x):=\phi(-x)$ for $x\in\R^2$. So it is not hard to check that each symmetry hyperplane of $\tilde{u}(x)$ and $\tilde{\phi}(x)$ contains the origin. Thus, repeating the above arguments for an arbitrary direction replacing the $x_1$-coordinate direction, we deduce that $\tilde{u}(x)$ and $\tilde{\phi}(x)$ are symmetric with respect to any hyperplane containing the origin, thus radially symmetric. Moreover, as a byproduct of the method, $\tilde{u}(x)$ and $\tilde{\phi}(x)$ are also strictly decreasing in the distance from the symmetry center. 	
\end{proof}

\section{Existence results for \eqref{Chs} by asymptotic approximation: proof of Theorem \ref{Thm:duozhjie}}\setcounter{equation}{0}\label{Sec_existence}

\noindent As we mentioned in the Introduction, the applicability of variational methods to the planar Choquard equation \eqref{Chs} is not straightforward. Indeed \eqref{Chs} has, at least formally, a variational structure related to the energy functional
\begin{multline}\label{eqn:log-fun}
	I(u):=\frac s2\int_{\R^2}\int_{\R^2}\frac{|u(x)-u(y)|^{\frac2s}}{|x-y|^4}\,{\rm d}x\,{\rm d}y
	+\frac s2\int_{\R^2}|u|^{\frac2s}\,{\rm d}x\\
	-\frac1{4\pi}\int_{\R^2}\int_{\R^2}\ln\frac1{|x-y|}F(u(y))F(u(x))\,{\rm d}x\,{\rm d}y\,.
\end{multline}
However, this energy functional is not well-defined on the natural Sobolev space $W^{s,\frac2s}(\R^2)$ because of 
the presence of the convolution term and the fact that the logarithm is unbounded both from below and from above. 
To overcome this difficulty, inspired by \cite{LRTZ,CDL} we will use an approximation technique as follows. Set
$$G_\alpha(x):=\frac{|x|^{-\alpha}-1}{\alpha},\qquad\alpha\in(0,1]\,,$$
$x\in\R^2$, and consider the modified approximating problem
\begin{equation}\label{eqn:2-spos}
    (-\Delta)^s_{\frac2s}u+|u|^{\frac2s-2}u=(G_\alpha(\cdot)\ast F(u))f(u)\quad\ \mbox{on}\ \ \R^2,
\end{equation}
with corresponding functional
\begin{equation*}
	\begin{split}
		I_\alpha(u):&=\frac s2\|u\|^{\frac2s}+\frac1{4\pi\alpha}\bigg[\int_{\R^2}F(u)\,{\rm d}x\bigg]^2-\frac1{4\pi\alpha}\int_{\R^2}\int_{\R^2}\frac1{|x-y|^{\alpha}}F(u(x))F(u(y))\,{\rm d}x\,{\rm d}y\\
		&=\frac s2\|u\|^{\frac2s}-\frac1{4\pi}\int_{\R^2}(G_\alpha(\cdot)\ast F(u))F(u)\,{\rm d}x\,.
	\end{split}
\end{equation*}
Unlike the original functional $I$, the power-type singularity in $G_\alpha$ can be handled by the Hardy-Littlewood-Sobolev inequality (Lemma \ref{Lem:HLS}), and it is standard to prove that $I_{\alpha}$ is well-defined and $C^1$ on $W^{s,\frac2s}(\R^2)$ with
$$
\aligned
	I'_\alpha(u)v=&\int_{\R^2}\int_{\R^2}\frac{|u(x)-u(y)|^{\frac2s-2}(u(x)-u(y))(v(x)-v(y))}{|x-y|^4}\,{\rm d}x\,{\rm d}y+\int_{\R^2}|u|^{\frac2s-2}uv\,{\rm d}x\\
	&-\frac1{2\pi}\int_{\R^2}(G_\alpha(\cdot)\ast F(u))f(u)v\,{\rm d}x
\endaligned
$$
for $u,v\in W^{s,\frac2s}(\R^2)$. In order to retrieve compactness in light of Lemma \ref{Lem:lions}, we will restrict  the functional $I_\alpha$ to the subspace $W_r^{s,\frac2s}(\R^2)$ of radially symmetric functions. Note that in the previous section we established symmetry of all strong positive solution for \eqref{Chs} belonging to $C^{1,1}_{\loc}(\R^2)$. Because of the singular behaviour of the $(s,\frac2s)$-fractional Laplacian, it seems difficult to prove regularity for weak solutions of \eqref{eqn:2-spos} belonging to $W^{s,\frac2s}(\R^2)$, so that the assumptions for the symmetry result in Theorem \ref{Thm:symmetric} can be fulfilled (see \cite{Brasco18}). However, since the problem is autonomous, it is natural to restrict $I_\alpha$ to the radial symmetric setting, by means of Palais' Principle of Symmetric Criticality (see \cite{Willem}), for which critical points of $I_\alpha$ restricted to $W_r^{s,\frac2s}(\R^2)$ are still unconstrained critical points for $I_\alpha$ on $W^{s,\frac2s}(\R^2)$.

\vskip0.2truecm
\noindent In the sequel we will use some elementary estimates which we collect in the next lemma. 
\begin{lemma}\label{Lem:G}
	Let $\alpha\in(0,1]$. Then,
	$$\frac{t^{-\alpha}-1}{\alpha}\geq\ln\frac1t\qquad\mbox{for all}\ \ t\in(0,1]\,.$$
	Moreover, for all $\nu>\alpha$ there exists $C_\nu>0$ such that
	$$
	\frac{t^{-\alpha}-1}{\alpha}\leq C_\nu t^{-\nu}\qquad\mbox{for all}\ \ t>0\,.
	$$
\end{lemma}

\noindent Let us show that for all $\alpha\in(0,1]$ the functional $I_\alpha$ enjoys a mountain-pass geometry.
\begin{lemma}\label{Lem:MP1}
	Let $\alpha\in(0,1]$ and assume \textcolor{red}{($f_1$)-($f_3$)}. Then, there exist constants $\rho,\eta>0$ and $e\in W^{s,\frac2s}_r(\R^2)$ such that:
	\begin{enumerate}
		\item[{\rm(i)}] $\|e\|>\rho$ and $I_\alpha(e)<0$;
		\item[{\rm(ii)}] $I_\alpha|_{S_\rho}\geq\eta>0$, where $S_\rho=\big\{u\in W_r^{s,\frac2s}(\R^2)\,|\,\|u\|=\rho\big\}$.
	\end{enumerate}
\end{lemma}
\begin{proof} (i)
	Take $e_0\in W_r^{s,\frac2s}(\R^2)$ such that $e_0(x)=1$ for $x\in B_\frac18(0)$, $e_0(x)=0$ for $x\in \R^2\setminus B_\frac14(0)$. For $t>0$ set
	\begin{equation}\label{eqn:AR-0}
		\Psi(t):=\frac12\bigg(\int_{\R^2}F(te_0)\,{\rm d}x\bigg)^2.
	\end{equation}
	By Remark \ref{Rmk_ass}-(ii) we infer that
	$$
	\frac{\Psi'(t)}{\Psi(t)}\geq\frac2{(s-\tau)t}\quad\ \text{for\,all}\ \ t>0
	$$
	and integrating over $[1,t]$ we find
	\begin{equation}\label{eqn:AR}
		\Psi(t)\geq\Psi(1)t^\frac2{s-\tau}=\frac12\left(\int_{\R^2}F(e_0)\,{\rm d}x\right)^2t^\frac2{s-\tau}.
	\end{equation}
	It follows from Lemma \ref{Lem:G}, \eqref{eqn:AR-0} and \eqref{eqn:AR} that
	\begin{equation*}
		\begin{split}
			I_{\alpha}(te_0)
			&\leq\frac s2t^\frac2s\|e_0\|^\frac2s-\frac1{4\pi}\int_{\R^2}\int_{\R^2}
			\frac{|x-y|^{-\alpha}-1}{\alpha}F\big(te_0(y)\big)F\big(te_0(x)\big)\,{\rm d}x\,{\rm d}y\\
			&=\frac s2t^\frac2s\|e_0\|^\frac2s-\frac1{4\pi}\int\int_{\{|x-y|\leq\frac12\}}\frac{|x-y|^{-\alpha}-1}{\alpha}F\big(te_0(y)\big)F\big(te_0(x)\big)\,{\rm d}x\,{\rm d}y\\
			&\leq\frac s2t^\frac2s\|e_0\|^\frac2s-\frac1{4\pi}\int\int_{{\{|x-y|\leq\frac12\}}}\ln\frac{1}{|x-y|}F\big(te_0(y)\big)F\big(te_0(x)\big)\,{\rm d}x\,{\rm d}y\\
			&\leq\frac s2t^\frac2s\|e_0\|^\frac2s-\frac{\ln2}{4\pi}\bigg(\int_{\R^2}F(te_0)\,{\rm d}x\bigg)^2\\
			&\leq\frac s2t^\frac2s\|e_0\|^\frac2s-\frac{\ln2}{4\pi}\left(\int_{\R^2}F(e_0)\,{\rm d}x\right)^2t^\frac2{s-\tau}.
		\end{split}
	\end{equation*}
	Therefore, one can find $t_0>0$ large enough such that $I_{\alpha}(t_0e_0)<0$.
	
	\medskip 
	
\noindent (ii) By Remark \ref{Rmk_ass}(i) there exists $C>0$ such that
	$$
	|F(u)|^{\frac43}\leq C\left(|u|^{\frac8{3s}}+|u|^{\frac8{3s}}\Phi_{2,s}\big(\tfrac43\alpha_*|u|^{\frac2{2-s}}\big)\right).
	$$
	Hence, by Hardy-Littlewood-Sobolev's inequality and Theorem \ref{Thm:m-theo}, we have
	\begin{equation}\label{eqn:guestimate}
		\aligned
			\int_{\R^2}\int_{\R^2}\frac1{|x-y|}&F(u(x))F(u(y))\,{\rm d}x\,{\rm d}y\leq C\bigg(\int_{\R^2}|F(u)|^{\frac43}\,{\rm d}x\bigg)^{\frac32}\\
			&\leq C\left(\|u\|_{\frac8{3s}}^{\frac4s}+\|u\|_{\frac8{3s}p'}^{\frac4s}\left(\int_{\R^2}\Phi_{2,s}\big(\tfrac43p\alpha_*\|u\|^{\frac2{2-s}}\left(\tfrac{|u|}{\|u\|}\right)^{\frac2{2-s}}\big){\rm d}x\right)^\frac3{2p}\right)
		\endaligned
	\end{equation}
	for some $p>1$ and $\frac1p+\frac1{p'}=1$. By Theorem \ref{Thm:m-theo}, we can find $p>1$ such that the last factor is bounded independently of $u$ provided $\|u\|<\big(\tfrac34\big)^{\frac{2-s}2}$. From (\ref{eqn:guestimate}) we deduce that
	\begin{equation*}
		\aligned
			I_\alpha(u)&=\frac s2\|u\|^{\frac2s}-\frac1{4\pi}\int_{\R^2}\int_{\R^2}\frac{|x-y|^{-\alpha}-1}{\alpha}F(u(y))F(u(x))\,{\rm d}x\,{\rm d}y\\
			&\geq\frac s2\|u\|^{\frac2s}-\frac1{4\pi}\int\int_{\{|x-y|\leq1\}}\frac{|x-y|^{-\alpha}-1}{\alpha}F(u(y))F(u(x))\,{\rm d}x\,{\rm d}y\\
			&\geq\frac s2\|u\|^{\frac2s}-\frac1{4\pi}\int\int_{\{|x-y|\leq1\}}\frac{F(u(y))F(u(x))}{|x-y|}\,{\rm d}x\,{\rm d}y\\
			&\geq\frac s2\|u\|^{\frac2s}-C(\|u\|_{\frac8{3s}}^{\frac4s}+\|u\|_{\frac8{3s}p'}^{\frac4s}).
		\endaligned
	\end{equation*}
So, let $\|u\|=\rho>0$ be sufficiently small, and recall the embedding $W_r^{s,\frac2s}(\R^2)\hookrightarrow L^t(\R^2)$ for all $t>\frac2s$, to obtain $\eta>0$ such that
	$I_{\alpha}(u)\geq \eta$ for any $\alpha \in (0,1]$.
\end{proof}

\noindent As a consequence of Lemma \ref{Lem:G}, the mountain pass level
$$
c_\alpha:=\inf\limits_{\gamma\in\Gamma}\max\limits_{t\in[0,1]}I_\alpha(\gamma(t)),
$$
where
$$
\Gamma:=\{\gamma\in C([0,1], W_r^{s,\frac2s}(\R^2)):\,\gamma(0)=0,\gamma(1)=e\}\,,
$$
turns out to be well defined. Moreover, from the Ekeland Variational Principle, the mountain pass geometry yields the existence of a Cerami sequence at level $c_\alpha$ for any fixed $\alpha\in(0,1]$, see e.g. \cite{Ekeland}. Namely, there exists $(u_n^\alpha)_n\subset W_r^{s,\frac2s}(\R^2)$ such that
$$
I_{\alpha}(u_n^\alpha)\to c_\alpha\quad\text{and}\quad(1+\|u_n^\alpha\|)I'_\alpha(u_n^\alpha)\to0\quad\text{in}\,\,(W^{s,\frac2s}(\R^2))'
$$
as $n\to+\infty$. For the sake of a lighter notation, we will simply use $u_n:=u_n^\alpha$.

\begin{remark}\label{rem:5.3}
	Observe from Lemma \ref{Lem:MP1} that there exist two constants $a,b>0$ independent of $\alpha$ such that $a<c_\al<b$.
\end{remark}

\noindent The next technical Lemma will be crucial in estimating the mountain pass level $c_\alpha$.
Let $R>0$ and $\bw\in C(\R^2)$ be such that
\begin{equation}\label{bar w}
	\bw(x):=\begin{cases}
		1\quad&\mbox{for}\ |x|\leq\frac R2,\\
		2-\frac2R|x|\quad&\mbox{for}\ |x|\in(\frac R2,R),\\
		0\quad&\mbox{for}\ |x|\geq R.
	\end{cases}
\end{equation}
\begin{lemma}\label{Lemma_bar w}
	For all $R>0$ and $s\in(0,1)$ we have $\bw\in W_r^{s,\frac2s}(\R^2)$ and
	$$\|\bw\|^\frac2s\leq\frac{\pi R^2}4\left(1+\frac{s(2+3s)}{(2+s)(1+s)}\right)+\frac{4\pi^2}{2-s}\left(\frac{31}{18}s+\frac{10}3\right).$$
\end{lemma}
\begin{proof}
	We have
	\begin{equation}\label{norm_bar w}
		\begin{split}
			\|\bw\|_\frac2s^\frac2s&=\frac{\pi R^2}4+2\pi\int_\frac R2^R\left|2-\frac2Rr\right|^\frac2sr\,\dd r=\frac{\pi R^2}4+2\pi\left(\frac2R\right)^\frac2s\int_0^\frac R2y^\frac2s(R-y)\,\dd y\\
			&=\frac{\pi R^2}4+2\pi\frac{R^2}2\left(\frac s{2+s}-\frac s{4+4s}\right)=\frac{\pi R^2}4\left(1+\frac{s(2+3s)}{(2+s)(1+s)}\right).
		\end{split}
	\end{equation}
	Let us now compute the seminorm of $\bw$. Since $\bw$ is radially symmetric, according to the equivalent formulations of the seminorm for radial functions established in \cite[Proposition 4.3]{Parini19}, we have
	\begin{equation*}
		\begin{split}
			[\bw]_{s,\frac2s}^\frac2s&=4\pi^2\int_0^{+\infty}\int_0^{+\infty}|\bw(r)-\bw(t)|^\frac2srt\frac{r^2+t^2}{|r^2-t^2|^3}\,{\rm d}r\,{\rm d}t\\
			&=8\pi^2\left(\int_0^\frac R2\int_\frac R2^R\ +\ \int_0^\frac R2\int_R^{+\infty}\ +\ \int_\frac R2^R\int_R^{+\infty}\right)+4\pi^2\int_\frac R2^R\int_\frac R2^R\\
			&=:8\pi^2\left(I_1+I_2+I_3\right)+4\pi^2I_4
		\end{split}
	\end{equation*}
and let us next estimate $I_i$, $i=1,\dots 4$, separately. Recalling that 
	\begin{equation}\label{derivative_PR}
		\frac\dd{\dd r}\left(\frac12\frac{r^2}{(t^2-r^2)^2}\right)=r\frac{r^2+t^2}{(t^2-r^2)^3}\,,
	\end{equation}
	we have
	\begin{equation*}
		\begin{split}
			I_1&=\int_\frac R2^R\left|\frac2Rt-1\right|^\frac2st\left(\int_0^\frac R2r\frac{r^2+t^2}{(t^2-r^2)^3}\,\dd r\right)\dd t=\frac12\int_\frac R2^R\left|\frac2Rt-1\right|^\frac2st\frac{R^2}4\frac{\dd t}{(t-\tfrac R2)^2(t+\tfrac R2)^2}\\
			&\leq\frac R8\left(\frac2R\right)^\frac2s\int_\frac R2^R\left|t-\frac R2\right|^{\frac2s-2}\,\dd t=\frac s{4(2-s)}\,.
		\end{split}
	\end{equation*}
	Again using \eqref{derivative_PR} we compute the second term $I_2$ as
	\begin{equation*}
		I_2=\int_0^\frac R2r\left(\int_R^{+\infty}t\frac{r^2+t^2}{(t^2-r^2)^3}\,\dd t\right)\dd r=\frac{R^2}2\int_0^\frac R2\frac r{(R^2-r^2)^2}\,\dd r=\frac13\,.
	\end{equation*}
	The third term $I_3$ can be estimated as
	\begin{equation*}
		\begin{split}
			I_3&=\int_\frac R2^R\left|\frac2Rr-1\right|^\frac2sr\left(\int_R^{+\infty}t\frac{r^2+t^2}{(t^2-r^2)^3}\dd t\right)\dd r\\
			&=\frac12\left(\frac2R\right)^\frac2s\int_\frac R2^R(R-r)^\frac2sr\frac{R^2}{(R-r)^2(R+r)^2}\,\dd r\\
			&\leq\frac29R\left(\frac2R\right)^\frac2s\int_\frac R2^R(R-r)^{\frac2s-2}\,\dd r=\frac49\frac s{2-s}\,.
		\end{split}
	\end{equation*}
Finally, integrating by parts we get
\begin{equation*}
		\begin{split}
			I_4&=\left(\frac2R\right)^\frac2s\int_\frac R2^Rr\left(\int_\frac R2^R |r-t|^{\frac2s}t\frac{r^2+t^2}{(t^2-r^2)^3}\,\dd t\right)\dd r\\
			&=\left(\frac2R\right)^\frac2s\int_\frac R2^Rr\bigg[|r-t|^{\frac2s}\cdot\frac12\frac{t^2}{(t^2-r^2)^2}\bigg]_{t=\frac{R}{2}}^{t=R}\,\dd r\\
			&+\left(\frac2R\right)^\frac2s\int_\frac R2^R\frac{r}{s}\int_\frac R2^r\bigg[(r-t)^{\frac2s-1}\cdot\frac{t^2}{(t^2-r^2)^2}\bigg]\dd t\,\dd r\\
			&-\left(\frac2R\right)^\frac2s\int_\frac R2^R\frac{r}{s}\int_r^R\bigg[(t-r)^{\frac2s-1}\cdot\frac{t^2}{(t^2-r^2)^2}\bigg]\dd t\,\dd r\\
			&=:A_1+A_2-A_3\,.
		\end{split}
	\end{equation*}
By inspection one has $A_3\geq0$. Moreover, 
\begin{equation*}
		\begin{split}
			A_1&\leq\left(\frac2R\right)^\frac2s\bigg[\int_\frac R2^R\frac{r}{2}|r-R|^{\frac2s-2}\frac{R^2}{(r+R)^2}\dd r-\int_\frac R2^R\frac{r}{2}|r-\tfrac{R}{2}|^{\frac2s-2}\frac{\frac14R^2}{(r+R/2)^2}\dd r\bigg]\\
			&\leq\left(\frac2R\right)^\frac2s\frac{R}{2}\int_\frac R2^R|r-R|^{\frac2s-2}\dd r=\frac{s}{2-s}
		\end{split}
	\end{equation*}
and
\begin{equation*}
\begin{split}
			A_2&=\left(\frac2R\right)^\frac2s\int_\frac R2^R\frac rs\bigg(\int_{\frac R2}^r(r-t)^{\frac2s-3}\cdot\frac{t^2}{(t+r)^2}\dd t\bigg)\dd r\\
			&\leq\left(\frac2R\right)^\frac2s\frac Rs\int_\frac R2^R(r-\tfrac R2)^{\frac2s-2}\,\dd r=\frac2{2-s}\,.
		\end{split}
\end{equation*}

\noindent Eventually we get  
	\begin{equation}\label{seminorm_bar w}
		[\bw]_{s,\frac2s}^\frac2s=4\pi^2(2I_1+2I_2+2I_3+I_4)\leq\frac{4\pi^2}{2-s}\left(\frac{31}{18}s+\frac{10}3\right).
	\end{equation}
	The desired estimate for the norm of $\bw$ is obtained combining \eqref{norm_bar w} with \eqref{seminorm_bar w}.
\end{proof}

\begin{lemma}\label{Lem:guji}
	Suppose (f$_1$)--(f$_5$) hold. Then $c_\alpha<\frac s4$ for all $\alpha\in(0,1]$.
\end{lemma}
\begin{proof}
	Let $R\in\big(0,\frac13\big]$ to be fixed later, and $\bw\in W_r^{s,\frac2s}(\R^2)$ be as in \eqref{bar w}. Define $w:=\|\bw\|^{-1}\bw$. As in the proof of Lemma \ref{Lem:MP1}, we deduce the existence of $T>0$ such that
	$$I_\alpha(Tw):=\max\limits_{t\geq0}I_\alpha(tw)\,,$$
	and therefore, since $I_\alpha$ is a functional of class $C^1$, one has $\frac{\rm d}{{\rm d}t}\big|_{t=T}I_\alpha(tw)=I'_\alpha(Tw)w=0$. Since by definition $c_\alpha\leq I_\alpha(Tw)$, we aim at proving that $I_\alpha(Tw)<\frac s4$. Suppose by contradiction that for all $\varepsilon>0$ there exists $\alpha\in(0,1]$ such that
	\begin{equation*}
		I_\alpha(Tw)\geq\frac s4\,,
	\end{equation*}
	then
	\begin{equation}\label{eqn:guji1}
		\aligned
		\frac s2\,T^{\frac2s}\geq\frac s4+\frac1{4\pi}\int_{\R^2}(G_\alpha(x)\ast F(Tw))F(Tw)\,{\rm d}x\,.
		\endaligned
	\end{equation}
	Observe that if $x,y\in B_R(0)$, with $R\leq\frac13$ then $G_\alpha(x-y)\geq\ln\frac1{|x-y|}\geq\ln\frac32>0$. Thus, since $w$ supported in $B_R(0)$, it follows from \eqref{eqn:guji1} that
	\begin{equation}\label{eqn:guji3-}
		T\geq\left(\frac12\right)^\frac s2.
	\end{equation}
	Moreover, $I_\alpha'(Tw)Tw=0$ yields
	\begin{equation}\label{eqn:guji2}
		T^{\frac2s}=\frac1{2\pi}\int_{\R^2}(G_\alpha(x)\ast F(Tw))f(Tw)Tw\,{\rm d}x\,.
	\end{equation}
	Let us estimate from below the right-hand side. We have
	\begin{equation}\label{eqn:guji4}
		\begin{split}
			&\int_{\R^2}(G_\alpha(x)\ast F(Tw))f(Tw)Tw\,{\rm d}x\\
			&\geq\int_{B_\frac R2\times B_\frac R2}G_\alpha(x-y)F(Tw(y))f(Tw(x))Tw(x)\,{\rm d}x\,{\rm d}y\\
			&\geq\int_{B_\frac R2\times B_\frac R2}\ln\frac1{|x-y|}F(Tw(y))f(Tw(x))Tw(x)\,{\rm d}x\,{\rm d}y\\
			&\geq\ln3\int_{|y|\leq\frac R2}F(Tw(y))\,\dd y\int_{|x|\leq\frac R2}f(Tw(x))Tw(x)\,{\rm d}x\\
			&\geq\ln3\bigg(\int_{|x|\leq\frac R2}\sqrt{F(Tw(x))f(Tw(x))Tw(x)}\,{\rm d}x\bigg)^2
		\end{split}
	\end{equation}
	by the H\"older inequality. Note that in $B_\frac R2(0)$ one has
	$$Tw=T\|\bar w\|^{-1}\geq\frac{\|\bar w\|^{-1}}{2^{\frac s2}}$$
	by \eqref{eqn:guji3-}. By the computations of the norm of $\bw$ carried out in Lemma \ref{Lemma_bar w}, one has $\|\bar w\|<\tfrac1{2^{\frac s2}T(s)}$ by taking $R=\frac13$ and
	\begin{equation}\label{def_T(s)}
		T(s):=\left(\frac\pi{18}\,\frac{4s^2+5s+2}{(2+s)(1+s)}+\frac{4\pi^2}9\,\frac{31s+60}{2-s}\right)^{-\frac s2}.
	\end{equation}
	Hence, by ($f_5$) we deduce from \eqref{eqn:guji2} and \eqref{eqn:guji4} that
	\begin{equation*}
		\begin{split}
			T^\frac2s&\geq\frac\beta{2\pi}\ln3\bigg(\int_{|x|\leq\frac R2}(Tw(x))^\frac{2+s}{2s}\,{\rm d}x\bigg)^2=\frac\beta{2\pi}\ln3|B_{\frac R2}(0)|^2\left(\frac T{\|\bw\|}\right)^{\frac2s+1}\\
			&=\beta\frac{\pi R^4\ln3}{32}\frac{T^{\frac2s+1}}{\|\bw\|^{\frac2s+1}}\,,
		\end{split}
	\end{equation*}
	which implies
	\begin{equation}\label{upperbound}
		T<\beta T\leq\frac{32\|\bw\|^{\frac2s+1}}{\pi R^4\ln3}\leq\frac{32}{\pi R^4\ln3}\left(\frac1{2^{\frac s2}T(s)}\right)^{\frac2s+1}=\frac{2^{4-\frac s2}3^4}{\pi\ln3}\left(\frac1{T(s)}\right)^{\frac2s+1}\ ,
	\end{equation}
	since $\beta>1$ and taking $R=\frac13$.
	Hence, combining \eqref{upperbound} and \eqref{eqn:guji3-}, we infer
\begin{equation*}
	 	\frac{2^{4-\frac s2}3^4}{\pi\ln3}\left(\frac1{T(s)}\right)^{\frac2s+1}\geq\beta T\geq\beta\left(\frac12\right)^\frac s2.
	 \end{equation*}
	This is in contradiction with ($f_5$) provided one takes
	\begin{equation}\label{def_beta}
		\beta>\beta_0:=\frac{2^4\,3^4}{\pi\ln3}\left(\frac1{T(s)}\right)^{\frac2s+1}.
	\end{equation}
\end{proof}

\noindent The estimate of the mountain-pass level obtained in Lemma \ref{Lem:guji} enables us to obtain the following compactness results 

\begin{lemma}\label{Lem:c-bounded}
	Assume that (f$_1$)--(f$_5$) hold. Let $(u_n)_n\subset W^{s,\frac2s}(\R^2)$ be a Cerami sequence of $I_{\alpha}$ at level $c_\alpha$, then $(u_n)_n$ is bounded in $W^{s,\frac2s}(\R^2)$ with
	\begin{equation}\label{eqn:bd6}
		\|u_n\|^\frac2s<\frac s\tau\,,
	\end{equation}
	as well as
	\begin{equation}\label{Lem:estimates_GFF_GFfu}
		\bigg|\int_{\R^2}[G_\alpha(x)\ast F(u_n)]F(u_n)\,{\rm d}x\bigg|<C\,,\quad\quad\bigg|\int_{\R^2}[G_\alpha(x)\ast F(u_n)]f(u_n)u_n\,{\rm d}x\bigg|<C\,.
	\end{equation}
\end{lemma}
\begin{proof}
	Since $(u_n)_n\subset W^{s,\frac2s}(\R^2)$ is a Cerami sequence, one has
	\begin{equation}\label{eqn:bd0}
		\frac s2\|u_n\|^\frac2s-\frac1{4\pi}\int_{\R^2}[G_\alpha(x)\ast F(u_n)]F(u_n)\,{\rm d}x\rightarrow c_\alpha
	\end{equation}
	and for all $v\in W^{s,\frac2s}(\R^2)$
	\begin{equation}\label{eqn:bd1}
		\aligned
			\int_{\R^2}\int_{\R^2}&\frac{|u_n(x)-u_n(y)|^{\frac{2-2s}s}(u_n(x)-u_n(y))(v(x)-v(y))}{|x-y|^4}\,{\rm d}x\,{\rm d}y\\	&+\int_{\R^2}|u_n|^{\frac2s-2}u_nv\,{\rm d}x-\frac1{2\pi}\int_{\R^2}[G_\alpha(x)\ast F(u_n)]f(u_n)v\,{\rm d}x=o_n(1)\|v\|\,.
		\endaligned
	\end{equation}
	Taking $v=u_n$, we get
	\begin{equation}\label{eqn:bd2}
		\int_{\R^2}\int_{\R^2}\frac{|u_n(x)-u_n(y)|^{\frac2s}}{|x-y|^4}\,{\rm d}x\,{\rm d}y+\int_{\R^2}|u_n|^{\frac2s}\,{\rm d}x -\frac1{2\pi}\int_{\R^2}[G_\alpha(x)\ast F(u_n)]f(u_n)u_n\,{\rm d}x=o_n(1)\,.
	\end{equation}
	In order to prove the boundedness of $(u_n)_n$, we introduce a suitable test function as follows
	$$
	\aligned
		v_n:=\left\{
		  \begin{array}{ll}
		    \frac{F(u_n)}{f(u_n)},&  u_n>0,\\
		    (s-\tau)u_n, & u_n\leq0,\\
		  \end{array}
		\right.
	\endaligned
	$$
	where $\tau$ is the positive constant appearing in $(f_3)$. It is easy to check that $|v_n|\leq C|u_n|$ since $F(t)\leq(s-\tau)f(t)t$ by Remark \ref{Rmk_ass}(ii) and $f(t)=0$ if and only if $t\leq0$. Furthermore,
	\begin{equation}\label{eqn:bd3}
		\begin{split}
			\int_{\R^2}\int_{\R^2}\frac{\left|v_n(x)-v_n(y)\right|^{\frac2s}}{|x-y|^4}&\,{\rm d}x\,{\rm d}y=\int\int_{u_n(x)>0 ; u_n(y)>0}\frac{\left|\frac{F\left(u_n(x)\right)}{f\left(u_n(x)\right)}-\frac{F\left(u_n(y)\right)}{f\left(u_n(y)\right)}\right|^{\frac2s}}{|x-y|^4}\,{\rm d}x\,{\rm d}y\\
			&+2\int\int_{u_n(x) \leq 0 ; u_n(y)>0} \frac{\left|(s-\tau) u_n(x)-\frac{F\left(u_n(y)\right)}{f\left(u_n(y)\right)}\right|^{\frac2s}}{|x-y|^4}\,{\rm d}x\,{\rm d}y\\
			&+(s-\tau)^\frac2s\int\int_{u_n(x)\leq 0; u_n(y) \leq 0}\frac{\left| u_n(x)-u_n(y)\right|^{\frac2s}}{|x-y|^4}\,{\rm d}x\,{\rm d}y\\
			&\leq\int_{\R^2}\int_{\R^2}\frac{\left(1-\frac{F\left(\xi_n(x, y)\right) f^{\prime}\left(\xi_n(x,y)\right)}{f^2\left(\xi_n(x,y)\right)}\right)^\frac2s\left|u_n(x)-u_n(y)\right|^{\frac2s}}{|x-y|^4}\,{\rm d}x\,{\rm d}y\\
			&+3(s-\tau)^{\frac2s}\int_{\R^2}\int_{\R^2} \frac{\left|u_n(x)-u_n(y)\right|^{\frac2s}}{|x-y|^4}\,{\rm d}x\,{\rm d}y\\
			&\leq\,C\int_{\R^2}\int_{\R^2} \frac{\left|u_n(x)-u_n(y)\right|^{\frac2s}}{|x-y|^4}\,{\rm d}x\,{\rm d}y\,,
		\end{split}
	\end{equation}
	by ($f_4$), where $\xi_n(x,y)\in\R$. Thus $v_n$ is well defined in $W^{s,\frac2s}(\R^2)$. Taking $v=v_n$ in (\ref{eqn:bd1}) and recalling that $f(t)=0$ and $F(t)=0$ for any $t\leq0$, we infer
	\begin{equation*}
	\aligned
		\int_{\R^2}\int_{\R^2}&\frac{|u_n(x)-u_n(y)|^{\frac2s-2}(u_n(x)-u_n(y))(v_n(x)-v_n(y))}{|x-y|^4}\,{\rm d}x\,{\rm d}y+(s-\tau)\int_{\{u_n<0\}} |u_n|^\frac2s\,{\rm d}x\\
		&+\int_{\R^2}|u_n|^{\frac2s-2}u_n\frac{F(u_n)}{f(u_n)}\,{\rm d}x-\frac1{2\pi}\int_{\R^2}[G_\alpha(x)\ast F(u_n)]F(u_n)\,{\rm d}x=o_n(1)\|u_n\|\,.
	\endaligned
	\end{equation*}
	Recalling (\ref{eqn:bd0}), this yields
	\begin{equation}\label{eqn:con9}
	\aligned
		\int_{\R^2}\int_{\R^2}&\frac{|u_n(x)-u_n(y)|^{\frac2s-2}(u_n(x)-u_n(y))(v_n(x)-v_n(y))}{|x-y|^4}\,{\rm d}x\,{\rm d}y\\
		&+\int_{\R^2}(s-\tau)|u_n|^\frac2s\,{\rm d}x+2c_\alpha-s\|u_n\|^\frac2s\geq o_n(1)\|u_n\|\,.
	\endaligned
	\end{equation}
	Similarly to (\ref{eqn:bd3}), we also obtain
	\begin{equation*}
	\begin{aligned}
	\int_{\R^2}&\int_{\R^2}\frac{\left|u_n(x)-u_n(y)\right|^{\frac2s-2}\left(u_n(x)-u_n(y)\right)\left(v_n(x)-v_n(y)\right)}{|x-y|^4}\,{\rm d}x\,{\rm d}y\\
	\leq &\int\int_{u_n(x)>0 ; u_n(y)>0} \frac{\left|u_n(x)-u_n(y)\right|^{\frac2s-2}\left(u_n(x)-u_n(y)\right)\left(\frac{F\left(u_n(x)\right)}{f\left(u_n(x)\right)}-\frac{F\left(u_n(y)\right)}{f\left(u_n(y)\right)}\right)}{|x-y|^4}\,{\rm d}x\,{\rm d}y\\
	&+2\int\int_{u_n(x) \leq 0 ; u_n(y)>0} \frac{\left|u_n(x)-u_n(y)\right|^{\frac2s-2}\left(u_n(x)-u_n(y)\right)\left((s-\tau) u_n(x)-\frac{F\left(u_n(y)\right)}{f\left(u_n(y)\right)}\right)}{|x-y|^4}\,{\rm d}x\,{\rm d}y\\
	&+\int\int_{u_n(x)\leq 0 ; u_n(y) \leq 0} \frac{(s-\tau)\left|u_n(x)-u_n(y)\right|^{\frac2s}}{|x-y|^4}\,{\rm d}x\,{\rm d}y\\
	\leq
	&(s-\tau)\int_{\R^2}\int_{\R^2}\frac{\left|u_n(x)-u_n(y)\right|^{\frac2s}}{|x-y|^4}\,{\rm d}x\,{\rm d}y\,.
	\end{aligned}
	\end{equation*}
	Combining this and (\ref{eqn:con9}) we get
	\begin{equation*}
		\tau\int_{\R^2}\int_{\R^2}\frac{|u_n(x)-u_n(y)|^{\frac2s}}{|x-y|^4}\,{\rm d}x\,{\rm d}y+\tau\int_{\R^2}|u_n|^{\frac2s}\,{\rm d}x\leq o_n(1)\|u_n\|+2c_\alpha\,.
	\end{equation*}
	As a consequence, we have
	\begin{equation*}
		\|u_n\|^\frac2s\leq\frac2\tau c_\alpha<\frac s\tau\,,
	\end{equation*}
	by \eqref{Lem:guji}, where we note that the constant on the right is independent of $n$ and $\alpha$. Finally, from (\ref{eqn:bd0}) and (\ref{eqn:bd2}) we immediately obtain \eqref{Lem:estimates_GFF_GFfu}.
\end{proof}

\begin{remark}\label{Rmk_nonneg_Cerami}
		Thanks to the uniform boundedness of Cerami sequences of Lemma \ref{Lem:c-bounded}, from now on we can always suppose that Cerami sequences at level $c_\alpha$ are nonnegative. Indeed, $u_n^-:=\min\{u_n,0\}\in W^{s,\frac2s}(\R^2)$ and thus, following the same argument used in Remark \ref{Rmk_positivity}, one has 
		\begin{equation*}
			\int_{\R^2}\int_{\R^2}\frac{|u_n^-(x)-u_n^-(y)|^\frac2s}{|x-y|^4}\,{\rm d}x\,{\rm d}y+\int_{\R^2}|u_n^-|^\frac2s\leq I'_\alpha(u_n)u_n^-\leq\|I'_\alpha(u_n)\|\|u_n\|=o_n(1)\,,
		\end{equation*}
		since $\|u_n\|\leq C$ by Lemma \ref{Lem:c-bounded}. This implies that $u_n^-\to0$ in $W^{s,\frac2s}(\R^2)$ as $n\to+\infty$ and therefore, setting $u_n^+:=\max\{u_n,0\}$, $(u_n^+)_n$ is a Cerami sequence of $I_\alpha$ at level $c_\alpha$, which we will denote simply by $(u_n)_n$.
\end{remark}

\begin{lemma}\label{Lem:integral-F}
Suppose ($f_1$)--($f_5$) hold and let $(u_n)_n$ be a Cerami sequence of $I_\alpha$ at level $c_\alpha$. Then there exists $C>0$ independent of $n$ and $\alpha$ such that
$$\int_{\R^2}f(u_n)u_n\,{\rm d}x\leq C \quad \text { and } \quad \int_{\R^2}F(u_n)^\kappa\,{\rm d}x\leq C$$
for any $\kappa\in\big[1,\tfrac1{\gamma_{s,\tau}}\big)$, where $\gamma_{s,\tau}\in(0,1)$ is a constant depending just on $s$ and $\tau$.
\end{lemma}
\noindent In the proof we will use the following 
\begin{lemma}
	Let $W\in[0,1)$ and $q>1$. Then
	\begin{equation}\label{estimate_1}
		(1+W)^q\leq 1+qW\max\{1,(1+W)^{q-1}\}
	\end{equation}
	and
	\begin{equation}\label{estimate_2}
		(1+W)^q-W-1\leq\frac{q!}{q-\lfloor q\rfloor}\frac W{1-W}\,.
	\end{equation}
\end{lemma}
\begin{proof}
	We have
	$$(1+W)^q=1+q\int_1^{1+W}t^{q-1}\dd t\leq 1+qW\max\{1,(1+W)^{q-1}\}\,.$$
	Next, let $N=\lfloor q\rfloor$. Hence, applying \eqref{estimate_1} iteratively, and from $\max\{1,(1+W)^{q-(N+1)}\}=1$, we get
	\begin{equation*}
		\begin{split}
			(1+W)^q-W-1&\leq1+qW(1+W)^{q-1}-W-1\\
			&\leq 1+qW+q(q-1)W^2+\cdots+q(q-1)\cdots(q-N)W^{N+1}-W-1\\
			&\leq(q-1)W+q(q-1)W^2+\cdots+q(q-1)\cdots(q-N)W^{N+1}\\
			&\leq q(q-1)\cdots(q-N+1)W(1+W+\cdots+W^N)\\
			&=\frac{q!}{q-\lfloor q\rfloor}W\frac{1-W^{N+1}}{1-W}\\
			&\leq\frac{q!}{q-\lfloor q\rfloor}\frac W{1-W}\,.
		\end{split}
	\end{equation*}
\end{proof}
\begin{proof}[Proof of Lemma \ref{Lem:integral-F}]
	Since $(u_n)_n$ is a Cerami sequence of $I_\alpha$ at level $c_\alpha$, by Lemma \ref{Lem:guji} we have
	\begin{equation}\label{eqn:integ2}
		\left\|u_n\right\|^{\frac2s}<\frac s\tau\,.
	\end{equation}
	Let us introduce the following auxiliary function
	$$
	G(t):=t-\frac1s\frac{F(t)}{f(t)}\quad\mbox{for}\ \,t \geq 0\,,
	$$
	and define $v_n:=G\left(u_n\right)$.
	Similarly to the proof of Lemma \ref{Lem:c-bounded}, one has $v_n \in W^{s,\frac2s}\left(\R^2\right)$.
	We first show that there exists $\gamma_{s,\tau}\in(0,1)$ depending just on $s$ and $\tau$, such that
	\begin{equation}\label{eqn:integ3-}
		\left\|v_n\right\|^{\frac2s}\leq\gamma_{s,\tau}<1
	\end{equation}
	for $n$ large enough. Indeed, recall by \eqref{eqn:bd0}-\eqref{eqn:bd1} that
	$$
	\frac1{2\pi}\int_{\R^2}\left[G_\alpha(x)\ast F\left(u_n\right)\right]F\left(u_n\right)\,{\rm d}x=s\left\|u_n\right\|^{\frac2s}-2c_\alpha+o_n(1)
	$$
	and
	$$
	\begin{aligned}
		&\int_{\R^2}\int_{\R^2} \frac{\left|u_n(x)-u_n(y)\right|^{\frac2s-2}\left(u_n(x)-u_n(y)\right)}{|x-y|^4}\left(\frac{F\left(u_n(x)\right)}{f\left(u_n(x)\right)}-\frac{F\left(u_n(y)\right)}{f\left(u_n(y)\right)}\right)\,{\rm d}x\,{\rm d}y\\
		&+\int_{\R^2}\left|u_n\right|^{\frac2s-2}u_n \frac{F\left(u_n\right)}{f\left(u_n\right)}\,{\rm d}x-\frac1{2\pi}\int_{\R^2}\left[G_\alpha(x)\ast F\left(u_n\right)\right] F\left(u_n\right)\,{\rm d}x=o_n(1),
	\end{aligned}
	$$
	which together imply
	\begin{equation}\label{eqn:integ3-bis}
		\begin{split}
			\frac1s&\int_{\R^2}\int_{\R^2} \frac{\left|u_n(x)-u_n(y)\right|^{\frac2s-2}\left(u_n(x)-u_n(y)\right)}{|x-y|^4}\left(\frac{F\left(u_n(x)\right)}{f\left(u_n(x)\right)}-\frac{F\left(u_n(y)\right)}{f\left(u_n(y)\right)}\right)\,{\rm d}x\,{\rm d}y\\
			&\quad+\frac1s\int_{\R^2}\left|u_n\right|^{\frac2s-2}u_n \frac{F\left(u_n\right)}{f\left(u_n\right)}\,{\rm d}x-\|u_n\|^\frac2s-\frac2sc_\alpha=o_n(1)\,.
		\end{split}
	\end{equation}
	Hence, one has 
	\begin{multline}\label{eqn:integ3}
		\|v_n\|^\frac2s=\int_{\R^2}\int_{\R^2}\frac{\left|G\left(u_n(x)\right)-G\left(u_n(y)\right)\right|^{\frac2s}}{|x-y|^4}\,{\rm d}x\,{\rm d}y+\int_{\R^2}\left|G\left(u_n\right)\right|^{\frac2s} \,{\rm d}x\\
		=\int_{\R^2}\int_{\R^2}\frac{\left|G\left(u_n(x)\right)-G\left(u_n(y)\right)\right|^{\frac2s}+\frac{1}{s}\left|u_n(x)-u_n(y)\right|^{\frac2s-2}\left(u_n(x)-u_n(y)\right)\left(\frac{F\left(u_n(x)\right)}{f\left(u_n(x)\right)}-\frac{F\left(u_n(y)\right)}{f\left(u_n(y)\right)}\right)}{|x-y|^4}\,{\rm d}y\,{\rm d}x\\
		+\int_{\R^2}|G(u_n)|^\frac2s\,{\rm d}x+\frac1s\int_{\R^2}|u_n|^{\frac2s-2}u_n\frac{F(u_n)}{f(u_n)}\,{\rm d}x+\frac2sc_\alpha-[u_n]_{s,\frac2s}^\frac2s-\|u_n\|_\frac2s^\frac2s+o_n(1)\,.
	\end{multline}
	Observe that for all $x,y\in\R^2$ by the mean value theorem one has
	$$\frac{F(u_n(x))}{f(u_n(x))}-\frac{F(u_n(y))}{f(u_n(y))}=\left(1-\frac{F(\theta(x,y))f'(\theta(x,y))}{f^2(\theta(x,y))}\right)\left(u_n(x)-u_n(y)\right),$$
	for some $\theta(x,y)\in\left(\min\left\{u_n(x),u_n(y)\right\},\max\left\{u_n(x),u_n(y)\right\}\right)$. Therefore, by the definition of $G$ we also get
	$$
		G\left(u_n(x)\right)-G\left(u_n(y)\right)=\left(1-\frac{1}{s}+\frac1s\frac{F(\theta(x, y))f^{\prime}(\theta(x, y))}{f^2(\theta(x, y))}\right)\left(u_n(x)-u_n(y)\right).
	$$
	Hence, we can estimate as follows:
	\begin{multline}\label{G_un_MVT}
		Z(u_n):=\\
		\int_{\R^2}\int_{\R^2}\frac{\left|G\left(u_n(x)\right)-G\left(u_n(y)\right)\right|^\frac2s+\frac1s\left|u_n(x)-u_n(y)\right|^{\frac2s-2}\left(u_n(x)-u_n(y)\right)\left(\frac{F\left(u_n(x)\right)}{f\left(u_n(x)\right)}-\frac{F\left(u_n(y)\right)}{f\left(u_n(y)\right)}\right)}{|x-y|^4}\,{\rm d}y\,{\rm d}x\\
		=\int_{\R^2}\int_{\R^2}\left[\left|1-\frac1s\left(1-\frac{F(\theta(x,y))f'(\theta(x,y))}{f^2(\theta(x,y))}\right)\right|^{\frac2s}+\frac1s\left(1-\frac{F(\theta(x,y))f'(\theta(x,y))}{f^2(\theta(x,y))}\right)\right]\,\cdot\\
		\cdot\, \frac{\left|u_n(x)-u_n(y)\right|^{\frac2s}}{|x-y|^4}\,{\rm d}x\,{\rm d}y\,.
	\end{multline}
	Let
	$$A:=\left\{(x,y)\in\R^2\times\R^2\,|\,\tfrac{Ff'}{f^2}(\theta(x,y))>1\right\}$$
	and split $\R^4=A\cup(\R^4\setminus A)$ and $Z(u_n)=Z_A(u_n)+Z_{\R^4\setminus A}(u_n)$ accordingly. On $\R^4\setminus A$ we have
	$$\left|1-\frac1s\left(1-\frac{F(\theta(x,y))f'(\theta(x,y))}{f^2(\theta(x,y))}\right)\right|^\frac2s\leq1-\frac1s\left(1-\frac{F(\theta(x,y))f'(\theta(x,y))}{f^2(\theta(x,y))}\right),$$
	and therefore
	\begin{equation}\label{I_R4-A}
		Z_{\R^4\setminus A}(u_n)\leq\int_{\R^4\setminus A}\frac{\left|u_n(x)-u_n(y)\right|^{\frac2s}}{|x-y|^4}\,{\rm d}x\,{\rm d}y\,.
	\end{equation}
	On the other hand, invoking the estimate \eqref{estimate_2} with $q=\frac2s$ and $W=W(x,y)=\frac1s\left(\frac{F(\theta(x,y))f'(\theta(x,y))}{f^2(\theta(x,y))}-1\right)>0$, we obtain
	\begin{equation*}
		\begin{split}
			Z_A(u_n)-\int_A\frac{\left|u_n(x)-u_n(y)\right|^{\frac2s}}{|x-y|^4}\,{\rm d}x\,{\rm d}y&=\int_A\left((1+W)^\frac2s-W-1\right)\frac{|u_n(x)-u_n(y)|^\frac2s}{|x-y|^4}\,{\rm d}x\,{\rm d}y\\
			&\leq\frac{\tfrac2s!}{\tfrac2s-\lfloor\tfrac2s\rfloor}\int_A\frac W{1-W}\frac{|u_n(x)-u_n(y)|^\frac2s}{|x-y|^4}\,{\rm d}x\,{\rm d}y\,.
		\end{split}
	\end{equation*}
	Requiring 
	\begin{equation}\label{def_mu_1}
		\frac{F(t)f'(t)}{f^2(t)}<1+\frac s2\quad\ \mbox{for all}\ \,t>0\,,
	\end{equation}
	one has $\|W\|_\infty<\tfrac12$. Hence, 
	\begin{equation}\label{I_A}
		Z_A(u_n)-\int_A\frac{\left|u_n(x)-u_n(y)\right|^{\frac2s}}{|x-y|^4}\,{\rm d}x\,{\rm d}y\leq\frac{2\tfrac2s!}{\tfrac2s-\lfloor\tfrac2s\rfloor}\|W\|_\infty\int_A\frac{|u_n(x)-u_n(y)|^\frac2s}{|x-y|^4}\,{\rm d}x\,{\rm d}y\,.
	\end{equation}
	Finally, from \eqref{I_R4-A}, \eqref{I_A} and \eqref{eqn:bd6}, we obtain
	\begin{equation*}\label{I_tot}
		\begin{split}
			Z(u_n)-[u_n]_{s,\frac2s}^\frac2s&\leq\frac{2\tfrac2s!}{\tfrac2s-\lfloor\tfrac2s\rfloor}\|W\|_\infty\int_A\frac{|u_n(x)-u_n(y)|^\frac2s}{|x-y|^4}\,{\rm d}x\,{\rm d}y\\
			&\leq\frac{2\tfrac2s!}{\tfrac2s-\lfloor\tfrac2s\rfloor}\|W\|_\infty\textcolor{red}{\frac s\tau}.
		\end{split}
	\end{equation*}
	We can now estimate the norm of $v_n$. From \eqref{eqn:integ3} we get
	\begin{equation*}
		\begin{split}
			\|v_n\|^\frac2s&=Z(u_n)-[u_n]_{s,\frac2s}^\frac2s+\int_{\R^2}|G(u_n)|^\frac2s\,{\rm d}x+\frac1s\int_{\R^2}|u_n|^{\frac2s-2}u_n\frac{F(u_n)}{f(u_n)}\,{\rm d}x\\
			&\quad+\frac2sc_\alpha-\|u_n\|_\frac2s^\frac2s+o_n(1)\\
			&\leq\frac{2\tfrac2s!}{\tfrac2s-\lfloor\tfrac2s\rfloor}\|W\|_\infty{\frac s\tau}+\int_{\R^2}|u_n|^\frac2s\left|1-\frac{F(u_n)}{s f(u_n)u_n}\right|^\frac2s{\rm d}x+\frac1s\int_{\R^2}\left|u_n\right|^{\frac2s-2}u_n\frac{F(u_n)}{f(u_n)}\,{\rm d}x\\
			&\quad+\frac2sc_\alpha-\|u_n\|_{\frac2s}^\frac2s+o_n(1)\\
			&\leq\frac{2\tfrac2s!}{\tfrac2s-\lfloor\tfrac2s\rfloor}\|W\|_\infty{\frac s\tau}+\int_{\R^2}|u_n|^\frac2s\left(1-\frac{F(u_n)}{sf(u_n)u_n}\right){\rm d}x+\frac1s\int_{\R^2}|u_n|^{\frac2s-2}u_n\frac{F(u_n)}{f(u_n)}\,{\rm d}x\\
			&\quad+\frac2sc_\alpha-\|u_n\|_\frac2s^\frac2s+o_n(1)\\
			&\leq\frac{\frac2s\tfrac2s!}{\tfrac2s-\lfloor\tfrac2s\rfloor}{\frac s\tau}\left(\left\|\frac{Ff'}{f^2}\right\|_\infty-1\right)+\frac12+o_n(1)\,,
		\end{split}
	\end{equation*}
	having used Remark \ref{Rmk_ass}(ii) and the estimate of $c_\alpha$ in Lemma \ref{Lem:guji}. In order to prove the claim \eqref{eqn:integ3-}, we need to check that the constant on the right-hand side is strictly less than $1$. This is equivalent to the following 
	\begin{equation}\label{def_mu_2}
		\left\|\frac{Ff'}{f^2}\right\|_\infty<1+\frac12\frac{\tfrac2s-\lfloor\tfrac2s\rfloor}{\frac2s\frac2s!}\,{\frac\tau s}.
	\end{equation}
	Defining now
	\begin{equation}\label{def_mu}
		\mu(s,\tau):=\min\left\{\frac s2,\,\frac12\frac{\tfrac2s-\lfloor\tfrac2s\rfloor}{\frac2s\frac2s!}\,{\frac\tau s}\right\},
	\end{equation}
	assumption ($f_3$) guarantees that both conditions \eqref{def_mu_1} and \eqref{def_mu_2} are satisfied.
	
	\noindent Next we aim at estimating the $L^1$ norm of a suitable power of $(F\left(u_n\right))_n$ by using \eqref{eqn:integ3-}. By $(f_4)$, for any $\varepsilon>0$ there exists $t_{\varepsilon}>0$ such that
	$$
	\frac{F(t)}{f(t)}\leq\varepsilon(t-t_\varepsilon)+\frac{F\left(t_{\varepsilon}\right)}{f\left(t_{\varepsilon}\right)}\quad\ \text{for all}\ \ t\geq t_{\varepsilon}\,.
	$$
	Hence, by $(f_3)$ one has
	\begin{equation}\label{eqn:integ4=}
		v_n=G\left(u_n\right)=\left(1-\frac{\varepsilon}{s}\right)\left(u_n-t_{\varepsilon}\right)+t_{\varepsilon}-\frac{1}{s} \frac{F\left(t_{\varepsilon}\right)}{f\left(t_{\varepsilon}\right)}\geq\left(1-\frac{\varepsilon}{s}\right)\left(u_n-t_{\varepsilon}\right),
	\end{equation}
	which implies that for all $x\in\R^2$
	\begin{equation}\label{eqn:integ4+}
		u_n(x)\leq t_\varepsilon+\frac{v_n(x)}{1-\bar{\varepsilon}}\,,
	\end{equation}
	where $\bar{\varepsilon}:=\frac{\varepsilon}{s}$.
	Hence, by $\left(f_1\right)$-$\left(f_2\right)$ we have that for any given $\varepsilon>0$, there exists $C_{\varepsilon}$ such that
	\begin{equation}\label{eqn:integ5}
		\begin{aligned}
			\int_{\R^2}F(u_n)^\kappa\,{\rm d}x&\leq\int_{|u_n|<t_{\varepsilon}}F(u_n)^\kappa\,{\rm d}x+\int_{|u_n|\geq t_{\varepsilon}}F(u_n)^\kappa\,{\rm d}x \\
			& \leq C_{\varepsilon}\left\|u_n\right\|^{\frac2s\kappa}+C_{\varepsilon}\int_{|u_n|\geq t_{\varepsilon}}\Phi_{2, s}\left(\alpha_*\kappa\left(t_{\varepsilon}+\frac{v_n}{1-\bar{\varepsilon}}\right)^{\frac{2}{2-s}}\right)\,{\rm d}x\,.
		\end{aligned}
	\end{equation}
	Moreover,
	$$
	\Phi_{2,s}\left(\alpha_*\kappa\left(t_{\varepsilon}+\frac{v_n}{1-\bar{\varepsilon}}\right)^{\frac2{2-s}}\right) \leq C_{\varepsilon} \Phi_{2, s}\left(\alpha_*\kappa(1+\bar{\varepsilon})\left(\frac{v_n}{1-\bar{\varepsilon}}\right)^{\frac{2}{2-s}}\right) .
	$$
	In view of \eqref{eqn:integ4=} and Remark \ref{Rmk_ass}(ii), $v_n\geq\frac\tau st_\varepsilon$ if $u_n\geq t_\varepsilon$, and then it follows from \eqref{eqn:integ5} that
	\begin{equation}\label{eqn:integ6}
	\int_{\R^2}F(u_n)^\kappa\,{\rm d}x\leq C_{\varepsilon}\|u_n\|^{\frac2s\kappa}+C_{\varepsilon}\int_{|u_n|\geq t_{\varepsilon}}
	\Phi_{2,s}\left(\alpha_*\kappa(1+\bar{\varepsilon})\left(\frac{v_n}{1-\bar{\varepsilon}}\right)^{\frac{2}{2-s}}\right){\rm d}x\,.
	\end{equation}
	Since $\left\|v_n\right\|^\frac2s\leq\gamma_{s,\tau}+\sigma<1$ for $n$ large enough and $\sigma>0$ small enough by \eqref{eqn:integ3-}, the following holds
	\begin{equation}\label{End_proof_Lemma_kappa}
		\frac{\kappa(1+\bar\varepsilon)}{(1-\bar\varepsilon)^{\frac2{2-s}}}\left\|v_n\right\|^{\frac2{2-s}}\leq\frac{\kappa(1+\bar\varepsilon)}{(1-\bar\varepsilon)^{\frac{2}{2-s}}}\left(\gamma_{s,\tau}+\sigma\right)<1
	\end{equation}
	for $\varepsilon>0$ small enough and $\kappa\in\big[1,\tfrac 1{\gamma_{s,\tau}}\big)$. As a consequence, from \eqref{eqn:integ6} we obtain
	$$
	\int_{\R^2}F(u_n)^\kappa\,{\rm d}x\leq C
	$$
	for some $C$ independent of $n$ and $\alpha$. Similarly, one can also prove that
	$$\int_{\R^2}f(u_n)u_n\,{\rm d}x\leq C$$
	for some $C$ independent of $n$ and $\alpha$.
\end{proof}

\medskip

\noindent Now we are in the position to prove the existence of a nontrivial critical point for the approximating functional $I_\alpha$, for $\alpha$ sufficiently small.

\begin{proposition}\label{Pro:compact}
	Assume that ($f_1$)--($f_5$) hold. For all $\alpha\in(0,1)$ sufficiently small there exists a positive $u_\alpha\in W^{s,\frac2s}_r(\R^2)$ such that $I_\alpha'(u_\alpha)=0$.
\end{proposition}
\begin{proof}
	By Lemma \ref{Lem:MP1} and Remark \ref{Rmk_nonneg_Cerami}, for all $\alpha\in(0,1)$ there exists a nonnegative Cerami sequence $(u_n^\alpha)_n\subset W^{s,\frac2s}_r(\R^2)$ of $I_\alpha$ at level $c_\alpha$, which is bounded in $W^{s,\frac2s}(\R^2)$ by Lemma \ref{Lem:c-bounded}. Hence there exists a nonnegative $u_\alpha\in W^{s,\frac2s}_r(\R^2)$ such that, up to a subsequence,
	\begin{equation}\label{eqn:fun3-0}
		\aligned
			u_n^\alpha\rightharpoonup u_\alpha\ \quad&\text{in}\,\,W^{s,\frac2s}(\R^2),\\
			u_n^\alpha\to u_\alpha\ \quad&\text{a.e. in}\,\, \R^2,\\
			u_n^\alpha\to u_\alpha\ \quad&\text{in}\,\, L^p(\R^2),\quad\mbox{for all}\ \ p\in\big(\tfrac2s,+\infty\big),
		\endaligned
	\end{equation}
	as $n\to+\infty$. Note that the compactness in the Lebesgue spaces is retrieved because we are restricting to the space of radially symmetric functions in $W^{s,\frac2s}(\R^2)$, see Lemma \ref{Lem:lions}. By the definition of $(u_n^\alpha)_n$ one has
	$$\|u_n^\alpha\|^2+\frac1{2\pi\alpha}\int_{\R^2}F(u_n^\alpha)\,{\rm d}x\int_{\R^2}f(u_n^\alpha)u_n^\alpha\,{\rm d}x=2c_\alpha+\frac1{2\pi\alpha}\int_{\R^2}\left(\frac1{|x|^\al}\ast F(u_n^\alpha)\right)f(u_n^\alpha)u_n^\alpha\,{\rm d}x+o_n(1)\,.$$
	Lemmas \ref{Lem:c-bounded} and \ref{Lem:integral-F} yield $C>0$ independent of $n$ and $\alpha$ such that
	\begin{equation}\label{eqn:fun2++}
		\int_{\R^2}\left(\frac1{|x|^\alpha}\ast F(u_n^\alpha)\right)f(u_n^\alpha(x))u_n^\alpha(x)\,{\rm d}x\leq C\,.
	\end{equation}
	Next we \textit{claim} that for any $\varphi\in C_0^\infty(\R^2)$
	\begin{equation}\label{eqn:fun3}
		\int_{\R^2}\left(\frac1{|x|^\al}\ast F(u_n^\alpha)\right)f(u_n^\alpha(x))\varphi(x)\,{\rm d}x\to\int_{\R^2}\left(\frac1{|x|^\alpha}\ast F(u_\alpha)\right)f(u_\alpha(x))\varphi(x)\,{\rm d}x\,,
	\end{equation}
	as $n\rightarrow\infty$. Indeed, first by Lemma \ref{Lem:integral-F}, for $x\in\R^2$ one has
	\begin{multline}\label{eqn:fun2+00}
			\int_{\R^2}\frac{F(u_n^\alpha(y))}{|x-y|^{\al}}\,{\rm d}y\leq\int_{|x-y|\leq 1}\frac{F(u_n^\alpha(y))}{|x-y|^{\al}}\,{\rm d}y+\int_{|x-y|>1}\frac{F(u_n^\alpha(y))}{|x-y|^{\al}}\,{\rm d}y\\
			\leq\bigg(\int_{|x-y|\leq 1}\frac{{\rm d}y}{|x-y|^{\frac32}}\bigg)^{\frac{2\alpha}3}\!\bigg(\int_{\R^2}F(u_n^\alpha(y))^{\frac3{3-2\alpha}}\,\dd y\bigg)^{\frac{3-2\alpha }3}\!+\int_{\R^2}F(u_n^\alpha(y))\,\dd y\leq C\,,
	\end{multline}
	provided $\al>0$ is small enough. Hence, in order to have \eqref{eqn:fun3}, it is sufficient to apply \cite[Lemma 2.1]{Figueiredo95} to the sequence 
	$$
	g(x,u_n^\alpha(x)):=\int_{\R^2}\frac1{|x-y|^\alpha}F(u_n^\alpha(y))\,{\rm d}y\,f(u_n^\alpha(x))
	$$
	restricted to any bounded domain $\Omega$. Indeed, $g(\cdot,u_n^\alpha)\in L^1(\Omega)$ by \eqref{eqn:fun2+00} and ($f_2$). Similarly, we can also prove
	\begin{equation}\label{eqn:fun3+}
		f(u_n^\alpha)\to f(u_\alpha)\quad\text{in}\,\,L_{loc}^1(\R^2)\,.
	\end{equation}
	Moreover, by \eqref{eqn:integ2}, \eqref{eqn:integ4+}, \eqref{eqn:integ3-}, and using the mean value theorem, we deduce
	\begin{equation}\label{eqn:fun2+0}
		\begin{split}
			\int_{\R^2}&|F(u_n^\alpha)-F(u_\alpha)|\,{\rm d}x
			=\int_{\R^2}|f(u_\alpha+\tau_n(x)(u_n^\alpha-u_\alpha))(u_n^\alpha-u_\alpha)|\,{\rm d}x\\
			&\leq {\varepsilon}\int_{\R^2}(|u_n^\alpha|+|u_\alpha|)^{\frac2s-1}|u_n^\alpha-u_\alpha|\,{\rm d}x\\
			&\quad+C_\varepsilon\int_{\R^2}\Phi_{2,s}\big(\alpha_*\big(u_\alpha+\tau_n(x)(u_n^\alpha-u_\alpha)\big)^{\frac{2}{2-s}}\big)|(u_n^\alpha+u_\alpha)(u_n^\alpha-u_\alpha)|\,{\rm d}x\\
			&\leq{\varepsilon\left(\|u_n^\alpha\|_\frac2s+\|u_\alpha\|_\frac2s\right)\|u_n^\alpha-u_\alpha\|_\frac2s}\\
			&\quad+ C_\varepsilon\bigg(\int_{\{u_n^\alpha>u_\alpha\}}\!\Phi_{2,s}\big(\alpha_*|u_n^\alpha|^{\frac2{2-s}}\big)|(u_n^\alpha+u_\alpha)(u_n^\alpha-u_\alpha)|\,{\rm d}x\\
			&\quad+\int_{\{u_n^\alpha\leq u_\alpha\}}\!\Phi_{2,s}\big(\alpha_*|u_\alpha|^{\frac2{2-s}}\big)|(u_n^\alpha+u_\alpha)(u_n^\alpha-u_\alpha)|\,{\rm d}x\bigg)\\
			&\leq {C\varepsilon}+C_\varepsilon\bigg(\int_{\{u_n^\alpha>u_\alpha\}}\Phi_{2,s}\big(\alpha_*	\Big|t_{\varepsilon}+\frac{v_n}{1-\bar{\varepsilon}}\Big|^{\frac{2}{2-s}}\big)|(u_n^\alpha+u_\alpha)(u_n^\alpha-u_\alpha)|\,{\rm d}x\\
			&\leq{C\varepsilon}+C\|u_n^\alpha\|_{\theta'\eta}\left(\int_{\R^2}\Phi_{2,s}\left(\alpha_*\theta'\eta'\Big|t_{\varepsilon}+\frac{v_n}{1-\bar{\varepsilon}}\Big|^{\frac{2}{2-s}}\right)\right)^\frac1{\theta'\eta'}\!\bigg(\int_{\R^2}|u_n^\alpha-u_\alpha|^\theta\,{\rm d}x\bigg)^\frac1\theta\\
			&\leq {C\varepsilon}+ C_\varepsilon\bigg(\int_{\R^2}|u_n^\alpha-u_\alpha|^\theta\,{\rm d}x\bigg)^{\frac{1}{\theta}}={C\varepsilon}+o_n(1)\,,
		\end{split}
	\end{equation}
	where $\tau_n(x)\in(0,1)$, $\bar\varepsilon=\tfrac\varepsilon s$, and for $\varepsilon$ small enough. Note that in the second-to-last inequality we have used H\"{o}lder's inequality with $\theta>\frac2s$ and $\theta'\eta>\frac2s$. Similarly, we can also infer
	\begin{equation}\label{eqn:fun2+01}
		\int_{\R^2}|f(u_n^\alpha)(u_n^\alpha-u_\alpha)|\,{\rm d}x\rightarrow0,\quad\text{as}\,\,n\rightarrow\infty\,.
	\end{equation}
	Consequently, combining \eqref{eqn:fun3-0}, \eqref{eqn:fun2+0}, \eqref{eqn:fun3} and \eqref{eqn:fun3+} we get $I'_\alpha(u_\alpha)=0$. Finally, we need to show that $u_\alpha$ is not trivial, by proving that $I_\alpha(u_\alpha)\neq0$. First, \eqref{eqn:fun2+00} and \eqref{eqn:fun2+01} imply
	$$
	\aligned
		\bigg|\int_{\R^2}\bigg[\bigg(\frac1{|x|^\al}&\ast F(u_n^\alpha)\bigg)f(u_n^\alpha(x))-\left(\frac1{|x|^\al}\ast F(u_\alpha)\right)f(u_\alpha(x))\bigg](u_n^\alpha(x)-u_\alpha(x))\,{\rm d}x\bigg|\\
		&\leq\bigg|\int_{\R^2}\bigg[\left(\frac{1}{|x|^\al}\ast F(u_n^\alpha)\right)-\left(\frac{1}{|x|^\al}\ast F(u_\alpha)\right)\bigg]f(u_n^\alpha(x))(u_n^\alpha(x)-u_\alpha(x))\,{\rm d}x\bigg|\\
		&+\bigg|\int_{\R^2}\left(\frac1{|x|^\al}\ast F(u_\alpha)\right)\big[f(u_n^\alpha(x))-f(u_\alpha(x))\big](u_n^\alpha(x)-u_\alpha(x))\,{\rm d}x\bigg|\\
		&=o_n(1)\,.
	\endaligned
	$$
	This, together with (\ref{eqn:fun2+0}) and (\ref{eqn:fun2+01}), implies that
	$$
	\aligned
		o_n(1)&=(I'_\alpha(u_n^\alpha)-I'_\al(u_\alpha))(u_n^\alpha-u_\alpha)\\
		&\geq C\|u_n^\alpha-u_\alpha\|^\frac2s+\frac1\alpha\bigg[\int_{\R^2}F(u_n^\alpha){\rm d}x\int_{\R^2}f(u_n^\alpha)u_n^\alpha\,{\rm d}x-\int_{\R^2}F(u_\alpha){\rm d}x\int_{\R^2}f(u_\alpha)u_\alpha\,{\rm d}x\bigg]\\
		&\quad+\frac1\al\bigg|\int_{\R^2}\bigg[\left(\frac1{|x|^\al}\ast F(u_n^\alpha)\right)f(u_n^\alpha(x))-\left(\frac1{|x|^\al}\ast F(u_\alpha)\right)f(u_\alpha(x))\bigg](u_n^\alpha(x)-u_\alpha(x))\,{\rm d}x\bigg|\\
		&=C\|u_n^\alpha-u_\alpha\|^\frac2s+o_n(1)\,.
	\endaligned
	$$
	Note that we used the monotonicity of the $(s,\frac2s)$-Laplacian, which relies on the well-known inequality
	$$(|a|^{p-2}a-|b|^{p-2}b)(a-b)\geq C_p|a-b|^p,$$
	which holds for all $p\geq2$ and $a,b\in\R^N$, see e.g. \cite{Simon}. As a consequence, $u_n^\alpha\to u_\alpha$ in $W_r^{s,\frac2s}(\R^2)$ as $n\to+\infty$, and therefore $I(u_\alpha)=c_\alpha\neq0$, which implies that $u_\alpha\neq0$. Hence, using the strong maximum principle for the $p-$fractional Laplacian \cite[Theorem 1.4]{Del17} we eventually obtain that $u_\alpha>0$ in $\R^2$.
\end{proof}

\medskip 

\noindent So far, we have proved that there exists $\bar\alpha\in(0,1)$ such that, for all $\alpha\in(0,\bar\alpha)$, the approximating functional $I_\alpha$ has a positive critical point in $W^{s,\frac2s}_r(\R^2)$. The last step consists in showing that actually the sequence $(u_\alpha)_\alpha$ converges, up to a subsequence, to a nontrivial critical point of the original functional $I$.

\noindent Thanks to the a-priori bound \eqref{eqn:bd6}, which is independent of $\alpha$, we have that $(u_\alpha)_\alpha$ is uniformly bounded in $W^{s,\frac2s}(\R^2)$ for $\alpha\in(0,\bar\alpha)$ with
\begin{equation}\label{eqn:re0}
	\|u_\alpha\|^\frac2s<\frac s\tau\,.
\end{equation}
Hence, up to a subsequence, there exists $u_0\in W_r^{s,\frac2s}(\R^2)$ such that
\begin{equation}\label{eqn:re1---}
	\aligned
		u_\alpha\rightharpoonup u_0\ \quad&\text{in}\,\, W^{s,\frac2s}(\R^2),\\
		u_\alpha\to u_0\ \quad&\text{a.e. in}\,\, \R^2,\\
		u_\alpha\to u_0\ \quad&\text{in}\,\, L^p(\R^2),\quad\mbox{for all}\ \ p\in\big(\tfrac2s,+\infty\big),
	\endaligned
\end{equation}
as $\alpha\to0^+$ by Lemma \ref{Lem:lions}.

\noindent Let us first prove some regularity results for $u_\alpha$. The strategy is similar to the one employed in Section \ref{Sec_equiv} and thus here we just highlight the differences.

\begin{lemma}\label{Lem:jianjin+}
	For any $\omega\in\big[1,\tfrac1{\gamma_{s,\tau}}\big)$ there exist $\bar\alpha\in(0,1)$ and a constant $C>0$, independent of $\alpha$, such that
	\begin{equation*}
		\int_{|x-y|\leq 1}\frac{F(u_\alpha(y))}{|x-y|^\frac{4(\omega-1)}{3\omega}}\,{\rm d}y\leq C\,,
	\end{equation*}
	and
	\begin{equation}\label{Lem:eq_infty}
		\int_{|x-y|\leq 1}\frac{F(u_\alpha(y))}{|x-y|^\frac{4(\omega-1)}{3\omega}}\,{\rm d}y\to 0\quad\ \ \mbox{as}\ \ |x|\to+\infty\,,
	\end{equation}
	uniformly for $\alpha\in(0,\bar\alpha)$.
\end{lemma}
\begin{proof}
	Let us define $v_\al:=G(u_\al)$, as in Lemma \ref{Lem:integral-F}. Using a similar argument, we deduce that $\sup_{\alpha\in(0,\bar\alpha)}\|v_\alpha\|<1$, and that for all $\omega\in\big[1,\tfrac1{\gamma_{s,\tau}}\big)$ there exists $C>0$ independent of $\alpha$ such that
	\begin{equation*}
		\int_{\R^2}F(u_\alpha)^\omega\,{\rm d}x<C\,.
	\end{equation*}
	Hence, by H\"older's inequality,
	\begin{equation}\label{F_omega}
		\int_{\{|x-y|\leq 1\}}\frac{F(u_\alpha(y))}{|x-y|^{\frac{4(\omega-1)}{3\omega}}}\,{\rm d}y\leq\bigg(\int_{\{|x-y|\leq1\}}\frac{{\rm d}y}{|x-y|^\frac43}\bigg)^{\frac{\omega-1}{\omega}}\bigg(\int_{\{|x-y|\leq1\}}F(u_\alpha)^\omega\,{\rm d}y\bigg)^\frac1\omega\leq C\,,
	\end{equation}
	where $C$ is independent of $\alpha$. Moreover, again as in Lemma \ref{Lem:integral-F}, for any $\varepsilon>0$ there exists $t_\varepsilon>0$ such that
	\begin{equation}\label{u_a v_a}
		u_\alpha(x)\leq t_\varepsilon+\frac{v_\alpha(x)}{1-\bar\varepsilon}\quad\text{for\,any}\,\,x\in\R^2,
	\end{equation}
	with $\bar\varepsilon:=\frac\varepsilon s$, which implies that there exists $C_\varepsilon>0$ such that
	$$
	u_\alpha^\frac2{2-s}\leq C_\varepsilon t_\varepsilon^\frac2{2-s}+(1+{\bar{\varepsilon}})^\frac2{2-s}\frac{v_\alpha^\frac2{2-s}}{(1-\bar{\varepsilon})^\frac2{2-s}}\,.
	$$
	Therefore, in order to show \eqref{Lem:eq_infty}, we can combine \eqref{F_omega} and \eqref{u_a v_a} and get
	$$
	\aligned
		&\int_{|x-y|\leq 1}\frac{F(u_\alpha(y))}{|x-y|^{\frac{4(\omega-1)}{3\omega}}}\,{\rm d}y\\
		&\leq C_\varepsilon\bigg(\int_{|x-y|\leq 1}|u_\alpha|^{\frac2s\omega}+|u_\alpha|^{\frac2s\omega}\Phi_{2,s}\bigg(\alpha_*\omega(1+\bar\varepsilon)^\frac2{2-s}
		 \frac{v_\al^\frac{2}{2-s}}{(1-\bar{\varepsilon})^\frac{2}{2-s}}\bigg){\rm d}y\bigg)^\frac1\omega\\
		&\leq C_\varepsilon\bigg(\int_{B_1(x)}|u_\alpha|^{\frac2s\omega}{\rm d}y+\big(\int_{B_1(x)}|u_\alpha|^{\frac2s\omega\zeta}{\rm d}y\big)^\frac1\zeta\bigg(\int_{\R^2}\Phi_{2,s}\bigg(\alpha_*\omega\zeta'(1+\bar\varepsilon)^\frac2{2-s}\frac{v_\alpha^\frac2{2-s}}{(1-\bar\varepsilon)^{\frac2{2-s}}}\bigg){\rm d}y\bigg)^\frac1{\zeta'}\bigg)^\frac1\omega\\
		&\leq C_\varepsilon\bigg(\int_{B_1(x)}|u_\alpha|^{\frac2s\omega}{\rm d}y+C_\varepsilon\big(\int_{B_1(x)}|u_\alpha|^{\frac2s\omega\zeta}{\rm d}y\big)^\frac1\zeta\bigg)^\frac1\omega,
	\endaligned
	$$
	where $\zeta'=\frac\zeta{\zeta-1}$, and in the last inequality we choose $\omega\zeta'$ sufficiently close to $1$ and we take advantage of the fact that $\|v_\alpha\|<1$. Thanks to \eqref{eqn:re1---} we obtain \eqref{Lem:eq_infty} and conclude the proof.
\end{proof}

\begin{lemma}\label{Lem:R2}
Let $\alpha\in(0,\bar\alpha)$ and $u_\alpha\in W^{s,\frac2s}(\R^2)$ be a weak solution of \eqref{eqn:2-spos}, then there exists $C>0$ independent of $\alpha$ such that $\|u_\alpha\|_\infty\leq C$. Furthermore, $u_\alpha\in C^{\nu}_{\loc}(\R^2)$ for some $\nu\in(0,1)$ and there exists $R>0$ such that for $|x|\geq R$,
$$
 |u_\alpha(x)|\leq\frac C{1+|x|^{\frac{7s}{2(2-s)}}}
$$
for some $C>0$ uniformly for $\alpha\in(0,\frac{4(\omega-1)}{3\omega})$.
\end{lemma}
\begin{proof}
	Since $u_\alpha$ is a nonnegative function, by (\ref{eqn:2-spos}) and Lemma \ref{Lem:G} we obtain
	\begin{equation}\label{eqn:re1_sez5}
	\aligned
		(-\Delta)^s_{\frac2s}{u_\alpha}+|u_\alpha|^{\frac2s-2}u_\alpha&\leq\left(\int_{|x-y|\leq 1}\frac{|x-y|^{-\alpha}-1}{\alpha}F(u_\alpha(y))\,{\rm d}y\right)f(u_\alpha(x))\\
		&\leq\int_{|x-y|\leq 1}\frac{F(u_\alpha(y))\,{\rm d}y}{|x-y|^{\frac{4(\omega-1)}{3\omega}}}\,f(u_\alpha(x))
	\endaligned
	\end{equation}
	weakly in $W^{s,\frac2s}(\R^2)$. For any $L>0$ and $\gamma>1$ consider the functions $\psi$ and $\Psi$ defined in \eqref{psi_Psi_K}. The proof follows the same line as in Lemma \ref{Lemma_infty}. Taking $\psi(u_\alpha)=u_\alpha u_{L,\alpha}^{\frac2s(\gamma-1)}$ as a test function in \eqref{eqn:re1_sez5}, we obtain
	\begin{multline}\label{eqn:re3_sez5}
	[\Psi(u_\al)]_{s,\frac2s}^{\frac2s}\\
    \leq\int_{\R^2}\int_{\R^2}\frac{|u_\alpha(x)-u_\alpha(y)|^{\frac2s-2}(u_\al(x)-u_\al(y))[u_\al(x)u_{L,\al}^{\frac{2(\gamma-1)}{s}}(x)-u_\al(y)u_{L,\al}^{\frac{2(\gamma-1)}{s}}(y)]}{|x-y|^4}\,{\rm d}y\,{\rm d}x\\
	\leq\int\int_{|x-y|\leq1}\frac{F(u_\alpha(y)){\rm d}y}{|x-y|^{\frac{4(\omega-1)}{3\omega}}}f(u_\alpha(x))u_\al u_{L,\al}^{\frac2s(\gamma-1)}\,{\rm d}x-\int_{\R^2}|u_\alpha|^{\frac2s}u_{L,\al}^{\frac2s(\gamma-1)}\,{\rm d}x\,.
	\end{multline}
	By ($f_1$)-($f_2$), for any $\epsilon>0$, there exists $C_\epsilon>0$ such that
	$$
	f(u_\al)\leq\epsilon|u_\al|^{\frac2s-1}+C_{\epsilon}|u_\al|^{\frac2s-1}\Phi_{2,s}(\alpha_*|u_\al|^{\frac2{2-s}}).
	$$
	So, for $\epsilon>0$ sufficiently small, using Young's inequality, Proposition \ref{Thm:m-theo}, Lemma \ref{Lem:jianjin+}, and (\ref{eqn:re3_sez5}), we get
	\begin{equation}\label{eqn:re4_sez5}
		{[\Psi(u_\al)]}^{\frac2s}_{s,\frac2s}\leq(C\epsilon-1)\|u_\al u_{L,\al}^{\gamma-1}\|_{\frac2s}^{\frac2s}+C_{\epsilon}\|u_\al u_{L,\al}^{\gamma-1}\|_{\frac{2p}{s}}^{\frac2s}\,,
	\end{equation}
	where $p$ is sufficiently large. Noting that all the constants are independent of $\alpha$, we can follow the argument of Lemma \ref{Lemma_infty} and obtain $\|u_\alpha\|_\infty\leq C$, with $C$ independent of $\alpha$. Finally, the decay estimate and H\"older's regularity can be obtained following step by step the proofs of Lemmas \ref{Lem:G-decay} and \ref{Lemma_alpha,0}, paying only attention to the fact that $G_\alpha$ substitutes the logarithmic Riesz kernel.
\end{proof}

\noindent We are finally in a position to prove Theorem \ref{Thm:duozhjie}, namely the existence of a weak solution of \eqref{Chs}. Next we prove that the sequence of solutions $(u_\alpha)_\alpha$ for the approximated problems \eqref{eqn:2-spos} has a nontrivial accumulation point $u_0\in W^{s,\frac2s}(\R^2)$ satisfying $I'(u_0)=0$.

\begin{proof}[Proof of Theorem \ref{Thm:duozhjie}]
	Let us divide the proof into two steps.
	
	\noindent \emph{Step 1.} We show that $u_0\in W^{s,\frac2s}_r(\R^2)$ satisfies $I'(u_0)=0$.\\
	For any $\varphi\in C_0^\infty(\R^2)$, we have
	\begin{multline}\label{eqn:Th1}
	I'_\alpha(u_\alpha)\varphi
	=\int_{\R^2}\int_{\R^2}\frac{|u_\al(x)-u_\alpha(y)|^{\frac2s-2}(u_\alpha(x)-u_\alpha(y))(\varphi(x)-\varphi(y))}{|x-y|^4}\,{\rm d}y\,{\rm d}x\\
	+\int_{\R^2}|u_\alpha|^{\frac2s-2}u_\alpha\varphi\,{\rm d}x-\frac1{2\pi}\int_{\R^2}\int_{\R^2}\frac{|x-y|^{-\alpha}-1}{\alpha}F(u_\alpha(y))\,{\rm d}y f(u_\alpha(x))\varphi\,{\rm d}x\,.
	\end{multline}
	The first two terms are easy to handle thanks to the convergence in \eqref{eqn:re1---}. Thus we focus on the third term and let us split $\R^2=\{|x-y|\leq1\}\cup\{|x-y|>1\}$. On the one hand, let $\omega\in\big(1,\frac1{\gamma_{s,\tau}}\big)$ be fixed, then by Lemmas \ref{Lem:G} and \ref{Lem:R2}, together with the continuity of $f$, for all $x,y\in\R^2$ such that $|x-y|\leq1$ we get
	\begin{equation}\label{eqn:Th2-}
		\aligned
			\left|\frac{|x-y|^{-\alpha}-1}{\alpha}F(u_\alpha(y))f(u_\alpha(x))\varphi(x)\right|&\leq	\frac{C_\omega}{|x-y|^\frac{4(\omega-1)}{3\omega}}F(u_\alpha(y))f(u_\alpha(x))|\varphi(x)|\\
			&\leq\frac C{|x-y|^\frac{4(\omega-1)}{3\omega}}|\varphi(x)|
		\endaligned
	\end{equation}
	for $\alpha<\tfrac{4(\omega-1)}{3\omega}$, which is readily integrable in $\{(x,y)\in\R^2\,|\,|x-y|\leq1\}$. Thus, recalling \eqref{eqn:re1---}, the Lebesgue dominated convergence theorem implies 
	\begin{equation}\label{eqn:Th2--}
		\aligned
			&\int\int_{|x-y|\leq1}\frac{|x-y|^{-\alpha}-1}\alpha F(u_\alpha(y))\,{\rm d}yf(u_\alpha(x))\varphi(x)\,{\rm d}x\\
			&\to\int\int_{|x-y|\leq1}\ln\frac1{|x-y|}\,F(u_0(y))\,{\rm d}yf(u_0(x))\varphi(x)\,{\rm d}x
		\endaligned
	\end{equation}
	as $\alpha\to0^+$. On the other hand, when $|x-y|\geq1$, there exists $\tau=\tau(|x-y|)\in(0,1)$ such that
	\begin{equation}\label{eqn:Th3}
	0\geq G_{\alpha}(x-y)=\frac{|x-y|^{-\alpha}-1}{\alpha}=-|x-y|^{-\tau\alpha}\ln|x-y|\,,
	\end{equation}
	applying the mean value theorem to the function $h_b(\alpha)=b^{-\alpha}$ with $b=|x-y|$. Since $\varphi$ has a compact 
support, it follows from Lemma \ref{Lem:R2}, ($f_1$), and the monotonicity of $F$, that
	\begin{equation}\label{eqn:Th4}
		\aligned
			\bigg|\frac{|x-y|^{-\alpha}-1}{\alpha}&F(u_\alpha(y))f(u_\alpha(x))\varphi(x)\bigg|\\
			&=|x-y|^{-\tau\alpha}\ln|x-y|F(u_\alpha(y))f(u_\alpha(x))|\varphi(x)| \\
			&\leq(|x|+|y|)F(u_\alpha(y))f(u_\alpha(x))|\varphi(x)|\\
			&\leq(1+|y|)\frac C{1+|y|^\frac7{2-s}}|\varphi(x)|\\
			&\leq\frac C{1+|y|^\frac52}|\varphi(x)|
		\endaligned
	\end{equation}
	which is readily integrable on $\{(x,y)\in\R^2\,|\,|x-y|\geq1\}$. Hence, again by the Lebesgue dominated convergence theorem, this implies
	\begin{equation}\label{eqn:Th5}
		\aligned
			&\int\int_{|x-y|\geq1}\frac{|x-y|^{-\alpha}-1}{\alpha}F(u_\alpha(y))\,{\rm d}yf(u_\alpha(x))\varphi(x)\,{\rm d}x\\
			&\to-\int\int_{|x-y|\geq1}\ln|x-y|F(u_0(y))\,{\rm d}yf(u_0(x))\varphi(x)\,{\rm d}x\,.
		\endaligned
	\end{equation}
	Finally, we need to prove the finiteness of 
	\begin{equation*}
		\int_{\R^2}\int_{\R^2}\ln|x-y|F(u_0(y))\,{\rm d}yF(u_0(x))\,{\rm d}x\,.
	\end{equation*}
By Fatou's lemma, we have
	\begin{equation}\label{eqn:Th6}
		\aligned
			&\bigg|\int_{\R^2}\int_{\R^2}\ln|x-y|F(u_0(y))\,{\rm d}y\,F(u_0(x))\,{\rm d}x\bigg|\\
			&\leq\liminf_{\alpha\to0}\bigg(\int\int_{|x-y|\leq1}G_{\alpha}(x-y)F(u_\alpha(y))\,{\rm d}y\,F(u_\alpha(x))\,{\rm d}x\\
			&-\int\int_{|x-y|\geq1}G_{\alpha}(x-y)F(u_\alpha(y))\,{\rm d}y\,F(u_\alpha(x))\,{\rm d}x\bigg).
		\endaligned
	\end{equation}
	Using Lemma \ref{Lem:G} with $\nu=\tfrac{4(\omega-1)}{3\omega}$, Lemma \ref{Lem:integral-F} with $\kappa=1$, and Lemma \ref{Lem:jianjin+}, we have
	\begin{equation}\label{eqn:Th7}
		\aligned
			\int\int_{|x-y|\leq1}G_\alpha&(x-y)F(u_\alpha(y))\,{\rm d}y\,F(u_\alpha(x))\,{\rm d}x\\
			&\leq C\int_{\R^2}F(u_\alpha(x))\left(\int_{|x-y|\leq1}\frac{F(u_\alpha(y))}{|x-y|^\frac{4(\omega-1)}{3\omega}}\,{\rm d}y\right)\,{\rm d}x\leq C
		\endaligned
	\end{equation}
	uniformly for $\alpha\in\big(0,\frac{4(\omega-1)}{3\omega}\big)$. So, by Remark \ref{rem:5.3} and (\ref{eqn:Th7}) we deduce
	\begin{equation}\label{eqn:Th8}
		\aligned
			&-\frac1{4\pi}\int\int_{|x-y|\geq1}G_\alpha(x-y)F(u_\alpha(y))\,{\rm d}y\,F(u_\alpha(x))\,{\rm d}x\\
			&\leq I_\alpha(u_\alpha)+\frac1{4\pi}\int\int_{|x-y|\leq1}G_\alpha(x-y)F(u_\alpha(y))\,{\rm d}y\,F(u_\alpha(x))\,{\rm d}x-\frac s2\|u_\alpha\|^\frac2s\leq C
		\endaligned
	\end{equation}
	uniformly for $\alpha$ sufficiently small. Combining \eqref{eqn:Th6}-\eqref{eqn:Th8}, we have
	\begin{equation}\label{eqn:Th9}
		\bigg|\int_{\R^2}\int_{\R^2}\ln|x-y|F(u_0(y))\,{\rm d}y\,F(u_0(x))\,{\rm d}x\bigg|<+\infty\,.
	\end{equation}
	Based on (\ref{eqn:Th2--}), (\ref{eqn:Th5}) and (\ref{eqn:Th9}), by taking the limit as $\alpha\to0^+$ in (\ref{eqn:Th1}), we have $I'(u_0)=0$ with $I(u_0)<+\infty$, that is, $u_0\in W^{s,\frac2s}_r(\R^2)$ solves equation \eqref{Chs}. \\
	
	\noindent \emph{Step 2.}
	To conclude the proof, we need to show that $u_0\neq0$ and that $u_\alpha\rightarrow u_0$ in $W^{s,\frac2s}(\R^2)$. Assume on the contrary that $u_\alpha\rightharpoonup 0$ in $W^{s,\frac2s}(\R^2)$, and so $u_\alpha\to0$ in $L^t(\R^2)$ for $t\in(\frac2s,+\infty)$. Similarly to (\ref{eqn:fun2+01}), we obtain $\int_{\R^2}f(u_\alpha)u_\alpha\,{\rm d}x=o_\alpha(1)$. Hence, by Lemma \ref{Lem:G} and Lemma \ref{Lem:jianjin+}, we have
	$$
	\aligned
		0=I'_\alpha(u_\alpha)u_\alpha&=\|u_\alpha\|^\frac2s-\frac1{2\pi}\int_{\R^2}\int_{\R^2}G_{\alpha}(x-y)F(u_\alpha(y))f(u_\alpha(x))u_\alpha(x)\,{\rm d}x\,{\rm d}y\\
		&\geq \|u_\alpha\|^\frac2s-\frac1{2\pi}\int\int_{\{|x-y|\leq1\}}G_{\alpha}(x-y)F(u_\alpha(y))f(u_\alpha(x))u_\alpha(x)\,{\rm d}x\,{\rm d}y\\
		&\geq \|u_\alpha\|^\frac2s-\frac1{2\pi}\int\int_{\{|x-y|\leq1\}}\frac{1}{|x-y|^{\frac{4(\omega-1)}{\omega}}}F(u_\alpha(y))f(u_\alpha(x))u_\alpha(x)\,{\rm d}x\,{\rm d}y\\
		&\geq \|u_\alpha\|^\frac2s-C\int_{\R^2}f(u_\alpha(x))u_\alpha(x)\,{\rm d}x\\
		&=\|u_\alpha\|^\frac2s+o_\alpha(1)\,,
	\endaligned
	$$
	which yields $u_\alpha\rightarrow0$ in $W^{s,\frac2s}(\R^2)$ as $\al\rightarrow0^+$. Then, using Remark \ref{rem:5.3} and \eqref{eqn:Th3}, we have
	$$
	\aligned
		a&\leq c_\alpha+o_\alpha(1)=I_\alpha(u_\alpha)\\
		&=\frac s2\|u_\alpha\|^\frac2s-\frac1{4\pi}\int_{\R^2}\int_{\R^2}G_{\alpha}(x-y)F(u_\al(y)) F(u_\al(x))\,{\rm d}x\,{\rm d}y\\
		&\leq-\frac1{4\pi}\int\int_{\{|x-y|\geq1\}}G_{\alpha}(x-y)F(u_\al(y))F(u_\al(x))\,{\rm d}x\,{\rm d}y+o_\al(1)\\
		&=-\frac1{4\pi}\int\int_{\{|x-y|\geq1\}}|x-y|^{-\tau(x,y)\alpha}\ln|x-y|F(u_\al(y))F(u_\al(x))\,{\rm d}x\,{\rm d}y+o_\al(1)\\
		&\leq\frac1{4\pi}\int_{\R^2}\int_{\R^2}\ln|x-y|F(u_\al(y))F(u_\al(x))\,{\rm d}x\,{\rm d}y+o_\al(1)\\
		&\leq \frac{\|F(u_\alpha)\|_1}{2\pi}\int_{\R^2}|x|F(u_\al(x))\,{\rm d}x+o_\al(1)\\
		&\leq C\bigg(R\int_{|x|\leq R}F(u_\al(x))\,{\rm d}x+\int_{|x|\geq R}|x|F(u_\al(x))\,{\rm d}x\bigg)+o_\al(1)\\
		&\leq C\bigg(\int_{|x|\leq R}|u_\al|^{\frac2s}\,{\rm d}x+\int_{|x|\geq R}\frac{C|x|}{1+|x|^\frac7{2-s}}\,{\rm d}x\bigg)+o_\al(1)\\
		&=o_\al(1)\,,
	\endaligned
	$$
	which yields a contradiction. Note that in the last inequality we have used Remark \ref{Rmk_ass}(i) with \eqref{eqn:re0} to estimate the first term, and ($f_1$) together with the monotonicity of $F$ and the decay given by Lemma \ref{Lem:R2}, to estimate the second term. Finally, similarly to \eqref{eqn:Th2--} and \eqref{eqn:Th5}, by Lemma \ref{Lem:R2} and the Lebesgue dominated convergence theorem, we have
	\begin{equation}\label{eqn:Th10}
		\aligned
			&\int_{\R^2}\int_{\R^2}\frac{|x-y|^{-\alpha}-1}{\alpha}F(u_\alpha(y))f(u_\alpha(x))u_\alpha(x)\,{\rm d}x\,{\rm d}y\\
			&\to-\int_{\R^2}\int_{\R^2}\ln|x-y|F(u_0(y))f(u_0(x))u_0(x)\,{\rm d}x\,{\rm d}y\,,
		\endaligned
	\end{equation}
	from which we conclude that $u_\alpha\rightarrow u_0$ in $W^{s,\frac2s}(\R^2)$ as $\alpha\to0^+$.
\end{proof}

\end{document}